\documentclass[journal]{IEEEtran} 

\usepackage{subfigure}
\usepackage{amsmath,amssymb}
\usepackage{url,booktabs}
\usepackage{graphicx}
\usepackage{balance}
\usepackage{color}
\usepackage{psfrag}
\usepackage{algorithm,algpseudocode}
\usepackage{rotating}
\usepackage{cite}
\usepackage[normalem]{ulem}
\usepackage{textcomp}
\usepackage[ugly]{units}

\graphicspath{{../../ImgFig/}}

\title{Comparing Edge Detection Methods based on Stochastic Entropies and Distances for PolSAR Imagery}

\author{Abra\~ao D.\ C.\ Nascimento,
        Michelle M.\ Horta,
        Alejandro C.\ Frery,~\IEEEmembership{Member}, 
        and
        \\ Renato J.\  Cintra,~\IEEEmembership{Senior Member}
\thanks{This work was supported by CNPq, Fapeal, FAPESP, and FACEPE, Brazil.}
\thanks{A.\ D.\ C.\ Nascimento is with the Departamento de Estat\'istica, Universidade Federal de Pernambuco, Cidade Universit\'aria, 50740-540, Recife, PE, Brazil, e-mail: abraao@de.ufpe.br}
\thanks{
M.\ M.\ Horta is with the Departamento de Computa\c c\~ao, Universidade Federal de S\~ao Carlos, Rod. Washington Lu\'{\i}s km 235, 13565-905, S\~ao Carlos, SP, Brazil, e-mail: michellemh@gmail.com}
\thanks{
A.\ C.\ Frery is with the Instituto de Computa\c c\~ao, Universidade Federal de Alagoas, BR 104 Norte km 97, 57072-970, Macei\'o, AL, Brazil, e-mail: acfrery@gmail.com}
\thanks{R.\ J.\ Cintra is with the Signal Processing Group, Departamento de Estat\'istica, Universidade Federal de Pernambuco, Cidade Universit\'aria, 50740-540, Recife, PE, Brazil, e-mail: rjdsc@ieee.org}
}

\begin{document}

\maketitle

\begin{abstract}
Polarimetric synthetic aperture radar (PolSAR) has achieved a prominent position as a remote imaging method. 
However, PolSAR images are contaminated by speckle noise due to the coherent illumination employed during the data acquisition.
This noise provides a granular aspect to the image,  making its processing and analysis (such as in edge detection)  hard tasks. 
This paper discusses seven methods for edge detection in multilook PolSAR images.
In all methods, the basic idea consists in detecting transition points in the finest possible strip of data which spans two regions. 
The edge is contoured using the transitions points and a B-spline curve. 
Four stochastic distances, two differences of entropies, and the maximum likelihood criterion were used under the scaled complex Wishart distribution; the first six stem from the $h$-$\phi$ class of measures.
The performance of the discussed detection methods was quantified and analyzed by the computational time and probability of correct edge detection, with respect to the number of looks, the backscatter matrix as a whole, the SPAN, the covariance an the spatial resolution.
The detection procedures were applied to three real PolSAR images. 
Results provide evidence that the methods based on the Bhattacharyya distance and the difference of Shannon entropies outperform the other techniques.
\end{abstract}

\section{Introduction}\label{sec:intro}

\IEEEPARstart{P}{olarimetric} synthetic aperture radar (PolSAR) has achieved a prominent position as a remote imaging method~\cite{LeePottier2009PolarimetricRadarImaging}.
SAR images are contaminated by speckle, a typical noise present in data acquired with coherent illumination subject to multipath interference.
Speckle noise introduces 
a granular aspect to the image, 
and its multiplicative nature makes SAR image analysis a challenging task~\cite{OliverandQuegan1998}. 
Most of the classical image processing methods assume additive noise, which is ineffective for processing SAR imagery~\cite{Gambinietal:IJRS:06}. 
Thus, PolSAR image analysis require specifically tailored signal processing techniques.

Among image processing techniques, edge detection occupies a central position. 
In simple terms, its purpose is to identify boundaries between regions of different structural characteristics~\cite{PolarimetricSegmentationBSplinesMSSP}.
In PolSAR systems, such features are usually represented by different scattering characteristics, 
which are effected by surface reflectance and speckle noise.

Several edge detection approaches have been proposed for SAR imagery, among them gradient-based methods~\cite{TouziLopesBousquet1988,OliverBlacknellWhite1996,FjortoftLopesMarthonCuberoCastan1998,XingyuHongjianKun2012}.
Such approach consists in using a sliding window to define a measure map that highlights the edges, resembling the edge strength map technique~\cite{TouziLopesBousquet1988}.
Then a thresholding step is applied to perform the sought edge detection.
In a different approach, Oliver~\textit{et al.}~\cite{OliverBlacknellWhite1996} proposed a maximum likelihood method aiming at two goals: (i)~detecting the presence of an edge within a window and (ii)~determining accurately the position of the edge.
Another line of research is the specific study of edge detection based on physical properties from urban areas.
As an example, 
Baselice and Ferraioli~\cite{BaseliceFerraioli2012} utilized Markov random fields to model jointly the amplitude and interferometric phase of two complex SAR images.

Statistical procedures, including active contour models, have also been applied to edge detection on SAR images~\cite{Horrit1999,GermainRefregier2001,Gambinietal:IJRS:06,GambiniandMejailandJacobo-BerllesandFrery,PolarimetricSegmentationBSplinesMSSP,GironFreryNeto2012}.
Giron~\textit{et al.}~\cite{GironFreryNeto2012} compare seven edge detectors following the general idea proposed by Gambini~\textit{et al.}~\cite{Gambinietal:IJRS:06}: finding transition points in strips of data.
This latter method was successfully employed in~\cite{GambiniandMejailandJacobo-BerllesandFrery} and~\cite{GironFreryNeto2012} for a comparative study using other strategies based on the same active contour approach.
In~\cite{PolarimetricSegmentationBSplinesMSSP}, Frery~\textit{et al.} extended it to multilook PolSAR data using the polarimetric $\mathcal{G}^H$ distribution as model.

In the present work, we aim at extending the proposal of Gambini~\textit{et al.}~\cite{Gambinietal:IJRS:06}. 
In its original idea, this method consists in forming strips of data around regions that span from the centroid of the candidate region to points located outside the region under consideration.
Resulting strips are then submitted to a screening phase, where each one is scanned looking for the point which defines two distinct regions under it based on maximization criterion, such as likelihood function.
The point that satisfies a decision rule is called \textit{transition point}.
A pre-processing step may be necessary in order to define the initial regions.
A post-processing step defines the edges of the regions by combining the transition points with B-spline curves.

For any type of approach, the edge detection problems can be related to three main aspects:
(i)~the procedure for detection, 
(ii)~the determination of the most accurate edge position, and 
(iii)~the specification of the window size (the window represents a square window or a strip of data).
The latter may influence all others aspects in such manner that smaller windows may not convey enough information to identify the presence of edges, while bigger windows may contain more than two edges. 
Then, the ideal window size is the one that provides only a single edge within the window~\cite{OliverBlacknellWhite1996}.
Following the approach discussed in~\cite{Gambinietal:IJRS:06}, we assume that there is one edge within the window provided by the initial selection. 

Our proposal improves or extends those procedures in three senses: 
(i)~the finest possible strips of data are used, namely the ones of one pixel width, and
(ii)~stochastic distances and difference of entropies are employed as objective functions to be maximized; not only the likelihood function as in~\cite{Gambinietal:IJRS:06}.
(iii)~the influence of the spatial image resolution is performed.

The Kullback-Leibler and Bhattacharyya distances have been used by Morio~\textit{et al.}~\cite{MorioRefregierGoudailFernandezDupuis2009} as a scalar contrast measure between different channels of polarimetric and interferometric synthetic aperture radar (PolInSAR) images.
These distances were calculated assuming that PolInSAR images follow the complex multidimensional Gaussian circular distribution.
The distances were then compared by their discrimination ability.
In the present paper, this work extended in two ways: 
(i)~the complex Wishart model is utilized instead of the Gaussian circular distribution, and 
(ii)~four distances, particular cases of the $h$-$\phi$ class of divergences, were considered as well as two difference of entropies.
Other information theory measures have been employed for change detection.
Erten~\textit{et al.}~\cite{Ertenetal2012} proposed a method based on mutual information for PolInSAR images assuming the scaled complex Wishart law.
As a result, they could verify that the new detector is more efficient than the one based on the
maximum likelihood ratio statistic~\cite{Conradsen2003}.

Therefore, the contribution of this paper is two fold. 
Firstly,
the paper presents a detailed discussion about information theory measures (using $h$-$\phi$ divergences and entropies) as criteria for the edge detection problem~\cite{Gambinietal:IJRS:06}.
Secondly, 
we compare their performance by edge accuracy and computational time, using contrast and spatial resolution as factors.
Three real multilook PolSAR images are employed to show the application of the proposal.
The best results were obtained using R\'enyi and Bhattacharya distances and difference of entropies.

This paper is structured as follows. 
Section~\ref{sec:1} presents the seven assessed measures. 
An analysis of the proposed techniques using real and simulated data is presented in Section~\ref{sec:2}.      
To that end, Monte Carlo experiments are performed in Section~\ref{application:first}  to assess required computational time and probability of correct detection.
Section~\ref{application:second} presents applications to real PolSAR images.
Finally, the results are summarized in Section~\ref{conclusion}.

\section{Edge detection approaches under scaled complex Wishart distribution}
\label{sec:1}

In this section, we present seven approaches for detecting edges in PolSAR images based on the methodology employed 
by~\cite{Gambinietal:IJRS:06,GambiniandMejailandJacobo-BerllesandFrery,PolarimetricSegmentationBSplinesMSSP,GironFreryNeto2012}: 
finding a transition point in a strip of data which is an estimator of the edge position.
In all cases, 
the estimator is the maximum of a given function.
The first function is a likelihood, 
which is a multivariate extension of the approach shown in~\cite{GambiniandMejailandJacobo-BerllesandFrery}.
Remaining functions are $h$-$\phi$ divergences and (difference of) entropies~\cite{FreryNascimentoCintraChileanJournalStatistics2011,salicruetal1994,salicruetal1993}. 

The edge detection methodologies considered in this paper 
operate in multiple stages:
(i)~identifying the initial centroid of the area of interest in a automatic, semiautomatic, or manual manner;
(ii)~casting rays from the centroid to the outside of the area;
(iii)~collecting data around the rays;
(iv)~detecting points in the strips of data which provide evidence of a change of properties, a transition;
(v)~defining the contour using a imputation method among the transition points, such as B-Splines~\cite{Gambinietal:IJRS:06}.
We bring contributions to stages~(iii) and~(iv).

\begin{figure}[hbt]
\centering
\includegraphics[width=1.03\linewidth]{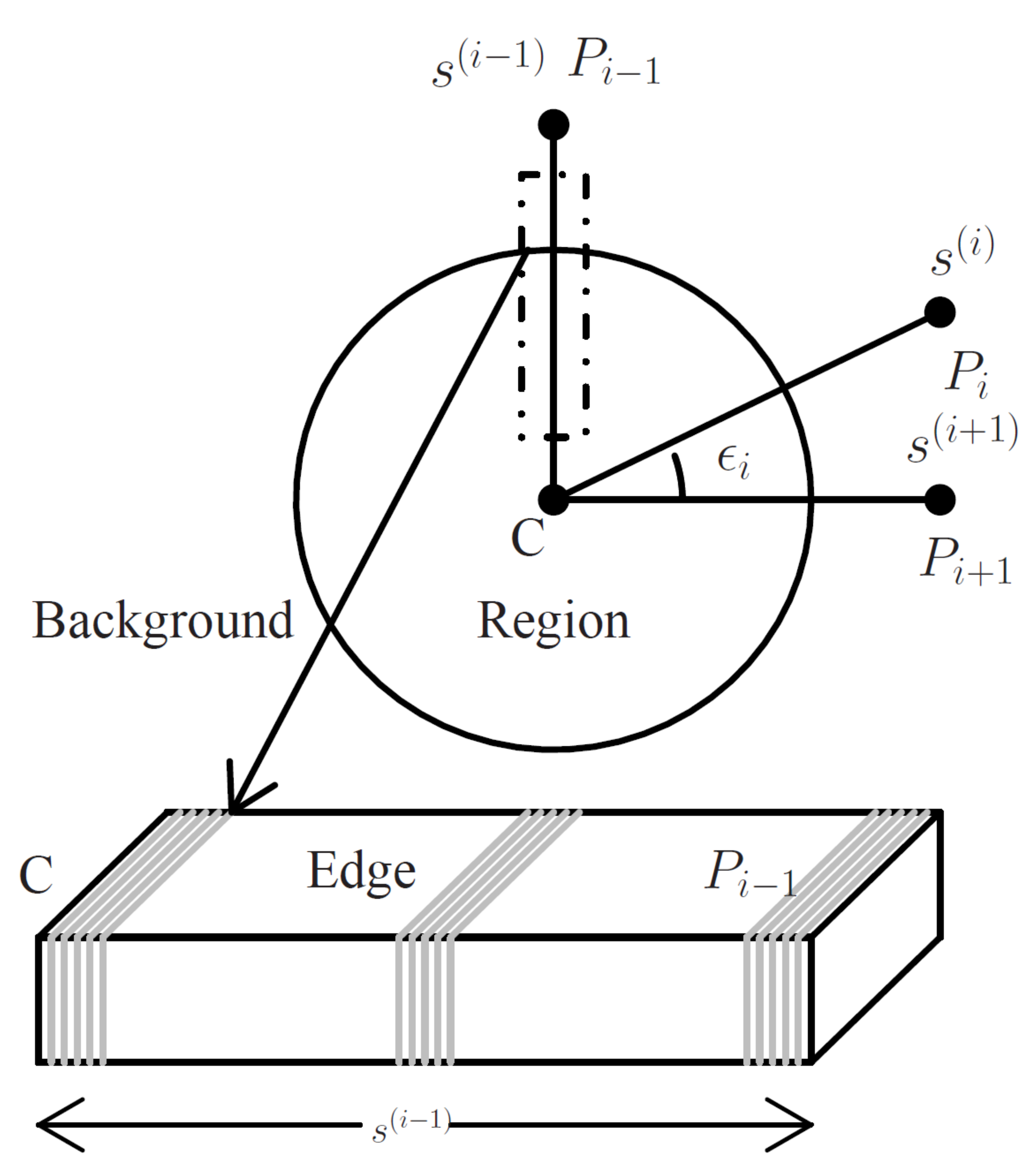}
\caption{
Edge detection on the polarimetric strip $s^{(i-1)}$ from the centroid C of a region to the control point $P_{i-1}$.
} 
\label{detectiongeometry}
\end{figure}

Initially,
let us admit 
a region $\mathcal{R}$ with centroid $C$.
Rays are traced from $C$ to control points $P_i$,
$i=1,2,\ldots,S$,
located outside $\mathcal R$;
the $S$ resulting rays can be represented by segments
in the form $\boldsymbol{s}^{(i)}=\overline{C P_i}$
and the angle between consecutive rays is 
$\epsilon_i=\angle(\boldsymbol{s}^{(i)},\boldsymbol{s}^{(i+1)})$.
Rays are converted into pixels using Bresenham's midpoint line algorithm~\cite{ComputerGraphicsMcConnell}.
This
representation provides the thinnest possible digital representation for a ray.
This contrast with the 20-pixel wide strips employed
in~\cite{Gambinietal:IJRS:06,GambiniandMejailandJacobo-BerllesandFrery,PolarimetricSegmentationBSplinesMSSP,GironFreryNeto2012}.
This setup is illustrated in Fig.~\ref{detectiongeometry}.

Data are assumed to follow a scaled complex Wishart distribution~\cite{FreitasFreryCorreia:Environmetrics:03}.
The scaled complex Wishart distribution
has density given by
\begin{equation}
 f_{\boldsymbol{Z}}(\boldsymbol{Z}';\boldsymbol{\Sigma},L) = \frac{L^{mL}|\boldsymbol{Z}'|^{L-m}}{|\boldsymbol{\Sigma}|^L \Gamma_m(L)} \exp\bigl[
-L\operatorname{tr}\bigl(\boldsymbol{\Sigma}^{-1}\boldsymbol{Z}'\bigr)\bigr],
\label{densitywish}
\end{equation}
where 
$\boldsymbol{Z}'$ is possible outcomes of $\boldsymbol{Z}$,
$\boldsymbol{\Sigma}$ represents the covariance matrix, 
$L$ is the number of looks, 
$m$ is the number of polarization channels, 
$\Gamma_m(L)=\pi^{m(m-1)/2}\prod_{i=0}^{L-1}\Gamma(L-i)$ is the multivariate gamma function, 
and $|\cdot|$ and $\operatorname{tr}(\cdot)$ are the determinant and the trace, respectively.
We refer to this distribution as ${\mathcal W}(\boldsymbol{\Sigma},L)$


The strip of data collected around the $i$th ray $\boldsymbol{s}^{(i)}$, 
$i=1,2,\ldots,S$, 
contains $N^{(i)}$ pixels.
Each pixel $k$ of a given strip $i$
is assumed
to be described by
the return matrix~$\boldsymbol{Z}^{(i)}_k$,
which is distributed according to the scaled complex Wishart law~\cite{PolarimetricSegmentationBSplinesMSSP}.
Denoting
the correct boundary position as $j^{(i)}$,
we have the following configuration:
\begin{equation}
\left\{
\begin{array}{rl}
\boldsymbol{Z}^{(i)}_k\sim\mathcal W(\boldsymbol{\Sigma}_A^{(i)},L_A^{(i)}),&\text{ for } k=1,\dots,j^{(i)}, \\
\boldsymbol{Z}^{(i)}_k\sim\mathcal W(\boldsymbol{\Sigma}_B^{(i)},L_B^{(i)}),&\text{ for } k=j^{(i)}+1,\dots,N^{(i)},
\end{array}
\right.\label{def1}
.
\end{equation}
In other words, each strip is formed by two patches of samples, and each patch obeys a complex Wishart law
with different parameter values.
In this work, 
we assume that the number of looks 
is estimated beforehand and 
remains
constant for the whole image.
Thus,
$L_A^{(i)}=L_B^{(i)}=L$ in every strip $i=1,2,\ldots,S$.
This supposition is realistic, 
since the number of looks is obtained 
from radar upon reception of 
backscattered pulses~\cite{FreitasFreryCorreia:Environmetrics:03}.

The main idea is to estimate the edge position $j^{(i)}$ 
on the data strip around ray $\boldsymbol{s}^{(i)}$ according to a specified decision rule.
Gray lines in Fig.~\ref{detectiongeometry} illustrate different configurations when the edge position is shifted. 

The model stated in~\eqref{def1} assumes that at most one transition occurs in any given ray.
Issues may arise if multiple transitions take place; this is discussed further in Section~\ref{application:second},
Fig.~\ref{otherana}.

In the following we present three different decision rules.
For brevity, we omit the superscript index $(i)$, since we focus our analysis on a single strip.

\subsection{Log-likelihood}

The likelihood function of the sample described by~\eqref{def1} is given by
\begin{align*}
L(j)
=
\prod_{k=1}^j 
f_{\boldsymbol{Z}}(\boldsymbol{Z}'_k;\boldsymbol{\Sigma}_A,L)   \prod_{k=j+1}^N  
f_{\boldsymbol{Z}}(\boldsymbol{Z}'_k;\boldsymbol{\Sigma}_B,L),
\end{align*}
where $\boldsymbol{Z}'_k$ is possible outcomes
of the random matrices described in~\eqref{def1}.
The log-likelihood is more convenient for maximization purposes:
\begin{equation}
\ell(j)=\log L(j)=\sum_{k=1}^{j}\log f_{\boldsymbol{Z}}(\boldsymbol{Z}'_k;\boldsymbol{\Sigma}_A,L) 
+\sum_{k=j+1}^N \log f_{\boldsymbol{Z}}(\boldsymbol{Z}'_k;\boldsymbol{\Sigma}_B,L).
\label{eq:loglik}
\end{equation}
The maximum likelihood estimator $\widehat{\jmath}_{\text{ML}}$ of the index 
at the location where the transition occurs is given by
\begin{equation}
\widehat{\jmath}_{\text{ML}}=\arg \displaystyle \max_j \ell(j).
\label{eq:jML}
\end{equation}

Using the model given by~\eqref{densitywish} in~\eqref{eq:jML}, 
with minor algebraic manipulation one obtains
\begin{align}
\ell(j)
=
&
N
\bigl[
-mL(1-\log L) - \log\Gamma_m(L)
\bigr]
+
L 
\bigl[
j\log|\widehat{\boldsymbol{\Sigma}}_A(j)|
\nonumber
\\
&\mbox{}
+
(N-j)
\log|\widehat{\boldsymbol{\Sigma}}_B(j)|
\bigr],
\label{loglik}
\end{align}
where 
$\widehat{\boldsymbol{\Sigma}}_A(j)$
and
$\widehat{\boldsymbol{\Sigma}}_B(j)$
are the maximum likelihood estimators for
$\Sigma_A$ and $\Sigma_B$, respectively,
with respect to edge position~$j$.
This estimator satisfies the following expression:
$$
\widehat{\boldsymbol{\Sigma}}_I(j)=
\left\{ \begin{array}{ll}
j^{-1}\sum_{k=1}^j \boldsymbol{Z}_k & \mbox{ if }I=A,\\
(N-j)^{-1}\sum_{k=j+1}^N \boldsymbol{Z}_k & \mbox{ if }I=B
\end{array} 
\right
.
$$


\subsection{Information Theory}

Information theory measures have received considerable attention as tools for contrast quantification~\cite{HypothesisTestingSpeckledDataStochasticDistances,InformationMeasuresInPerspective}.
In this section, we describe novel methodologies for detecting PolSAR boundaries based on such measures.
We adopt the following notation.

Let $\boldsymbol{X}$ and $\boldsymbol{Y}$ 
be two random matrices obeying the scaled complex Wishart distribution
with same known numbers of looks.
The associated probability densities 
are $f_{\boldsymbol{X}}(\boldsymbol{Z}';\boldsymbol{\theta}_1)$ and $f_{\boldsymbol{Y}}(\boldsymbol{Z}';\boldsymbol{\theta}_2)$, respectively, 
where $\boldsymbol{\theta}_1=\operatorname{vec}(\boldsymbol{\Sigma}_A)$ and $\boldsymbol{\theta}_2=\operatorname{vec}(\boldsymbol{\Sigma}_B)$ are their parameter vectors and $\operatorname{vec}(\cdot)$ is the vectorization operator. 
Both densities are defined over the set of Hermitian  positive definite matrices $\mathcal A$.

\subsubsection{Stochastic distances}

\begin{itemize}
\item [(i)] Kullback-Leibler:
\begin{align*} 
d_{\text{KL}}(\boldsymbol{X},\boldsymbol{Y}) = \frac 12\int_{\mathcal A}&[f_{\boldsymbol{X}}(\boldsymbol{Z}';\boldsymbol{\theta}_1)-f_{\boldsymbol{Y}}(\boldsymbol{Z}';\boldsymbol{\theta}_2)]\\
&\times\log\frac{f_{\boldsymbol{X}}(\boldsymbol{Z}';\boldsymbol{\theta}_1)}{f_{\boldsymbol{Y}}(\boldsymbol{Z}';\boldsymbol{\theta}_2)}\,\,\mathrm{d}\boldsymbol{Z}'.
\end{align*}
\item [(ii)] R\'{e}nyi of order $0<\beta<1$:
\begin{align*}
d_{\text{RD-}\beta}(\boldsymbol{X},\boldsymbol{Y}) =&\frac1{\beta-1}\\
\times\log\Bigl[&\frac{ \int_{\mathcal A} {f_{\boldsymbol{X}}(\boldsymbol{Z}';\boldsymbol{\theta}_1)}^{\beta}{f_{\boldsymbol{Y}}(\boldsymbol{Z}';\boldsymbol{\theta}_2)}^{1-\beta}\,\,\mathrm{d}\boldsymbol{Z}'}2\\ 
&+\frac{\int_{\mathcal A} {f_{\boldsymbol{X}}(\boldsymbol{Z}';\boldsymbol{\theta}_1)}^{1-\beta}{f_{\boldsymbol{Y}}(\boldsymbol{Z}';\boldsymbol{\theta}_2)}^{\beta}\,\,\mathrm{d}\boldsymbol{Z}'}2\Bigr].
\end{align*}
\item [(iii)] Bhattacharyya:
$$
d_{\text{BA}}(\boldsymbol{X},\boldsymbol{Y}) = -\log\int_{\mathcal A}\sqrt{f_{\boldsymbol{X}}(\boldsymbol{Z}';\boldsymbol{\theta}_1)f_{\boldsymbol{Y}}(\boldsymbol{Z}';\boldsymbol{\theta}_2)}\,\,\mathrm{d}\boldsymbol{Z}'.
$$
\item [(iv)] Hellinger:
$$
d_{\text{H}}(\boldsymbol{X},\boldsymbol{Y}) = 1-\int_{\mathcal A}\sqrt{f_{\boldsymbol{X}}(\boldsymbol{Z}';\boldsymbol{\theta}_1)f_{\boldsymbol{Y}}(\boldsymbol{Z}';\boldsymbol{\theta}_2)}\,\,\mathrm{d}\boldsymbol{Z}'.
$$
\end{itemize}
In the above expressions, the differential element is given by
$$
\mathrm{d}\boldsymbol{Z}'=\prod_{i=1}^m\mathrm{d}Z'_{ii}\displaystyle\prod^m_{\underbrace{i,j=1}_{i<j}}\mathrm{d}\Re\{Z'_{ij}\} \mathrm{d}\Im\{Z'_{ij}\},
$$ 
where $Z'_{ij}$ is the $(i,j)$th entry of matrix $\boldsymbol{Z}'$, and $\Re\{\cdot\}$ and $\Im \{\cdot\}$ return the real and imaginary parts of their arguments, respectively.
Provided that the densities characterize the same distribution and possibly differ only in the parameter,
we use
$(\boldsymbol{\theta}_1,\boldsymbol{\theta}_2)$
instead of
$(\boldsymbol{X},\boldsymbol{Y})$
as arguments
of the considered distances.
%
%
Frery~\textit{et al.}~\cite{FreryCintraNascimento2013} derived 
analytical expressions 
for the above distances
considering 
Wishart distributions for the general case 
$\boldsymbol{\Sigma}_A\neq\boldsymbol{\Sigma}_B$ and $L_A\neq L_B$.
Since we assume that the equivalent number of looks is the same, 
the distances can be expressed according to:
\begin{itemize}
\item[(i)] Kullback-Leibler:
$$ 
d_\text{KL}(\boldsymbol{\theta}_1,\boldsymbol{\theta}_2)=
L\bigg[\frac{\operatorname{tr}(\boldsymbol{\Sigma}_1^{-1}\boldsymbol{\Sigma}_2+\boldsymbol{\Sigma}_2^{-1}\boldsymbol{\Sigma}_1)}{2}-
m\bigg].  
$$
\item[(ii)] R\'enyi:
\begin{align*}
d_{\text{RD}-\beta}(\boldsymbol{\theta}_1,\boldsymbol{\theta}_2)=
&\frac{\log 2}{1-\beta}+\frac{1}{\beta-1}\log\Big\{ \\
&\Bigl[\frac{|(\beta \boldsymbol{\Sigma}_1^{-1}+(1-\beta)\boldsymbol{\Sigma}_2^{-1})^{-1}|}{|\boldsymbol{\Sigma}_1|^{\beta}|\boldsymbol{\Sigma}_2|^{1-\beta}}\Bigr]^L \\
+
&\Bigl[\frac{|(\beta \boldsymbol{\Sigma}_2^{-1}+(1-\beta)\boldsymbol{\Sigma}_1^{-1})^{-1}|}{|\boldsymbol{\Sigma}_1|^{(1-\beta)}|\boldsymbol{\Sigma}_2|^{\beta}}\Bigr]^L\Big\}.
\end{align*}
\item[(iii)] Bhattacharyya
\begin{align*}
d_\text{BA}(\boldsymbol{\theta}_1,\boldsymbol{\theta}_2)=&
L\bigg[\frac{\log|\boldsymbol{\Sigma}_1|+\log|\boldsymbol{\Sigma}_2|}{2} \\
&\mbox{}-\log\bigg|\bigg(\frac{\boldsymbol{\Sigma}_1^{-1}+\boldsymbol{\Sigma}_2^{-1}}{2}\bigg)^{-1}\bigg|\bigg].
\end{align*}
\item[(iv)] Hellinger
$$
\hskip-5em
d_\text{H}(\boldsymbol{\theta}_1,\boldsymbol{\theta}_2) =1-\Bigg[\frac{\bigl|2^{-1}({\boldsymbol{\Sigma}_1^{-1}+\boldsymbol{\Sigma}_2^{-1}})^{-1}\bigr|}{\sqrt{|\boldsymbol{\Sigma}_1||\boldsymbol{\Sigma}_2|}}\Bigg]^L. 
$$
\end{itemize}
For $L=1$,
the Kullback-Leibler and Bhattacharyya distances become
similar to 
the expressions derived by Morio~\textit{et al.}~\cite{MorioRefregierGoudailFernandezDupuis2009}.

Above expressions are modified according to~\cite{salicruetal1994}
in order to obtain test statistics with known asymptotic properties.
Thus
test statistics based on these distances can be readily derived 
for the null hypothesis
$\mathcal H_0:\boldsymbol{\theta}_1=\boldsymbol{\theta}_2$:
\begin{equation}
S_{\mathcal D}\big(\widehat{\boldsymbol{\theta}_1}{(j)},\widehat{\boldsymbol{\theta}_2}{(N-j)}\big) = \frac{2 j (N-j)v_{\mathcal D}}{N}
d_{\mathcal D}\big(\widehat{\boldsymbol{\theta}_1}{(j)},\widehat{\boldsymbol{\theta}_2}{(N-j)}\big),
\label{eq:TestStatistic}
\end{equation}
where ${\mathcal D}\in\{\text{KL, RD-}\beta\text{, BA, H}\}$, $v_{\mathcal D} \in\{ 1,\beta^{-1}, 4, 4\}$, respectively, $\widehat{\boldsymbol{\theta}_1}{(j)}=\operatorname{vec}(\widehat{\boldsymbol{\Sigma}}_A(j))$ and $\widehat{\boldsymbol{\theta}_2}{(N-j)}=\operatorname{vec}(\widehat{\boldsymbol{\Sigma}}_B(N-j))$ are estimators for $\boldsymbol{\theta}_1$ and $\boldsymbol{\theta}_2$ using random samples of sizes $j$ and $N-j$ as well as , respectively.

Let $p$ be the number of elements of the parameter vectors, i.e. $p=m^2$, as discussed in~\cite{EstimationEquivalentNumberLooksSAR}.
As derived by Frery~\textit{et al.}~\cite{FreryCintraNascimento2013}, 
the test statistics given in~\eqref{eq:TestStatistic} are asymptotically $\chi^2$-distributed 
with $p$ degrees of freedom.
This result holds true as long as
the
maximum likelihood estimators are used with large samples~\cite{salicruetal1994}.

Thus, 
four novel detectors can be proposed 
for finding edges on PolSAR imagery 
by seeking for the point that 
maximizes the test statistics between the two models, i.e.,
\begin{align*}
\widehat{\jmath}_{\mathcal D}=\arg \displaystyle \max_j 
\underbrace{S_{\mathcal D}\big(\widehat{\boldsymbol{\theta}_1}{(j)},\widehat{\boldsymbol{\theta}_2}{(N-j)}\big)}_{\equiv \,\,S_{\mathcal D}(j)}
=\arg \displaystyle \max_j S_{\mathcal D}(j),
\end{align*}
where $\mathcal D=\{\text{KL},\text{RD-}\beta,\text{BA},\text{H}\}$.

\subsubsection{Stochastic entropies}

The concepts of ``information'' and ``entropy'' received its fundamental mathematical treatment in the context of data communications by Shannon in 1948~\cite{Shannon1948}.
Thenceforth, the proposition and the application of these measures have become active research fields in several areas.
In particular, Frery~\textit{et al.}~\cite{EntropyBasedStatisticalAnalysisPolSAR} derived entropies and 
associated test statistics for the scaled complex Wishart model.

As suggested by Pardo~\textit{et al.}~\cite{Pardo1997},
we consider the R\'enyi of order $0<\beta<1$  
and Shannon entropies, 
denoted by
$H_{\text{R-}\beta}$ and $H_{\text{S}}$,
respectively.
A discrimination statistics based on stochastic entropies can be defined as
\begin{align}
S_{\mathcal M}(\widehat{\boldsymbol{\theta}_1}{(j)},\widehat{\boldsymbol{\theta}_2}{(N-j)})=&\frac{j\big(H_{\mathcal M}(\widehat{\boldsymbol{\theta}_1}{(j)})-\overline{v}\big)^2}{\sigma_{\mathcal M}^2(\widehat{\boldsymbol{\theta}_1}{(j)})}\nonumber \\
&+\frac{(N-j)\big(H_{\mathcal M}(\widehat{\boldsymbol{\theta}_2}{(N-j)})-\overline{v}\big)^2}{\sigma_{\mathcal M}^2(\widehat{\boldsymbol{\theta}_2}{(N-j)})},
\label{eq:EntropyTest}
\end{align}
where 
$H_{\mathcal M}\in\{H_{\text{S}},H_{\text{R-}\beta}\}$, 
$\boldsymbol{\theta}_i=[\theta_{i1}\;\theta_{i2}\;\cdots\;\theta_{ip}]^\top$,
$i=1,2$,
$\theta_{ik}$ is the $k$th element of vector $\boldsymbol{\theta}_i$, $k=1,2,\ldots,p$,
$\sigma_{\mathcal M}^2(\boldsymbol{\theta}_i)
=
\boldsymbol{\delta}_i^* 
\mathcal K(\boldsymbol{\theta}_i)^{-1}
\boldsymbol{\delta}_i$, 
the superscript
${}^\top$ is the transposition operator,
${}^*$ represents the conjugate transposition,
$\mathcal K(\boldsymbol{\theta}_i)=\operatorname{E}\{-\partial^2 \log f_{\boldsymbol{Z}}(\boldsymbol{Z};\boldsymbol{\theta}_i)/\partial \boldsymbol{\theta}_i^2\}$ is the Fisher information  matrix, 
$\boldsymbol{\delta}_i=[\delta_{i1}\;\delta_{i2}\;\cdots\;\delta_{ip}]^\top$ such that $\delta_{ik}=\partial H_{\mathcal M}(\boldsymbol{\theta}_i)/\partial \theta_{ik}$, $p$ is the size of vector $\boldsymbol{\theta}_i$, 
and
\begin{align*}
\overline{v}=&\Biggl[
\frac j{\sigma_{\mathcal M}^2(\widehat{\boldsymbol{\theta}_1}{(j)})}+ \frac {(N-j)}{\sigma_{\mathcal M}^2(\widehat{\boldsymbol{\theta}_2}{(N-j)})}
\Biggr]^{-1}\\
&\times\Biggl[ 
\frac{j H_{\mathcal M}(\widehat{\boldsymbol{\theta}_1}{(j)})}{\sigma_{\mathcal M}^2(\widehat{\boldsymbol{\theta}_1}{(j)})}+\frac{(N-j) H_{
\mathcal M}(\widehat{\boldsymbol{\theta}_2}{(N-j)})}{\sigma_{\mathcal M}^2(\widehat{\boldsymbol{\theta}_2}{(N-j)})}
\Biggr].
\end{align*}
Frery~\textit{et al.}~\cite{EntropyBasedStatisticalAnalysisPolSAR} showed that the expression given in~\eqref{eq:EntropyTest} follows asymptotically a $\chi^2_1$ distribution, and it can be used to test whether two samples from the scaled Wishart distribution possess the same entropy, i.e., ${\mathcal H}_0:H_{\mathcal M}({\boldsymbol{\theta}_1})=H_{\mathcal M}({\boldsymbol{\theta}_2})$.
Additionally, the Shannon and R\'enyi of order $\beta$ entropies are given by, respectively:
\begin{align*}
H_{\text{S}}(\boldsymbol{X}) \equiv H_{\text{S}}(\boldsymbol{\theta}_i)=&-\int_{\boldsymbol{\mathcal{A}}} f_{\boldsymbol{X}}(\boldsymbol{Z}';\boldsymbol{\theta}_i)\log f_{\boldsymbol{X}}(\boldsymbol{Z}';\boldsymbol{\theta}_i) \mathrm{d}\boldsymbol{Z}'\\
=&\operatorname{E}\{-\log f_{\boldsymbol{X}}(\boldsymbol{Z};\boldsymbol{\theta}_i)\}
\end{align*}
and
\begin{align*}
H_{\text{R-}\beta}(\boldsymbol{X}) \equiv H_{\text{R-}\beta}(\boldsymbol{\theta}_i)=&(1-\beta)^{-1}\log\int_{\boldsymbol{\mathcal{A}}} f_{\boldsymbol{X}}^\beta(\boldsymbol{Z}';\boldsymbol{\theta}_i)\mathrm{d}\boldsymbol{Z}'\\
=&(1-\beta)^{-1}\log\operatorname{E}\bigl\{f_{\boldsymbol{X}}^{\beta-1}(\boldsymbol{Z};\boldsymbol{\theta}_i)\bigr\},
\end{align*}
for $i=1,2$.
Applying these expressions to the density given in~\eqref{densitywish} the following entropies and variances are obtained:
\begin{itemize}
\item Shannon:
\begin{align*}
H_\text{S}(\boldsymbol{\Sigma},L)=&\frac{m(m-1)}{2}\log \pi- m^2 \log L +m \log |\boldsymbol{\Sigma}| \\
&+mL+(m-L)\psi_m^{(0)}(L) 
+\displaystyle \sum_{k=0}^{m-1}\log \Gamma(L-k)
\end{align*}
and
\begin{align*}
\sigma_{\text{S}}^2=&\frac{\bigl[(m-L)\psi_m^{(1)}(L)+m-\frac{m^2}{L}\bigr]^2}{\psi_m^{(1)}(L) - \frac{m}{L}}\\
&+\frac{m^2}{L} \operatorname{vec}\bigl(\boldsymbol{\Sigma}^{-1}\bigr)^* \bigl( \boldsymbol{\Sigma} \otimes \boldsymbol{\Sigma}\bigr)  \operatorname{vec}\bigl(\boldsymbol{\Sigma}^{-1}\bigr),
\end{align*}
where $\otimes$ is the Kronecker product, ``$\operatorname{vec}$'' is the operator which vectorizes its argument, and $\psi_m^{(v)}$ is the $v$th-order multivariate polygamma function given by~\cite{EstimationEquivalentNumberLooksSAR}
$$
\psi_m^{(v)}(x)=\sum_{i=0}^{m-1} \psi^{(v)}(x-i),
$$
and the ordinary polygamma function is~\cite{EstimationEquivalentNumberLooksSAR}
$$
\psi^{(v)}(x)=\frac{\partial^{v+1} \log\Gamma(x)}{\partial x^{v+1}},
$$
for $v\geq 0$.

\item R\'enyi of order $0<\beta<1$:
Denoting $q=L+(1-\beta)(m-L)$, the R\'enyi entropy is expressed by
\begin{align*}
H_{\text{R-}\beta}(\boldsymbol{\Sigma}_I,L)=& \frac{m(m-1)}{2}\log\pi- m^2 \log L +m \log |\boldsymbol{\Sigma}_I| \\
&+\frac{\sum_{i=0}^{m-1}\bigr[\log\Gamma(q-i)-\beta\log\Gamma(L-i)\bigl]}{1-\beta}\\
&-\frac{mq\log\beta}{1-\beta}
\end{align*}
and
\begin{align*}
\sigma_{\text{R-}\beta}^2=&
\frac{\Big\{\frac{\beta}{1-\beta} \bigl[\psi_m^{(0)}(q)-\psi_m^{(0)}(L)\bigr]-\frac{m\beta\log(\beta)}{1-\beta}-\frac{m^2}{L}\Big\}^2}{\psi_m^{(1)}(L) - \frac{m}{L}}\\
&+\frac{m^2}{L} \operatorname{vec}\bigl(\boldsymbol{\Sigma}^{-1}\bigr)^* \bigl(\boldsymbol{\Sigma}_I \otimes \boldsymbol{\Sigma}\bigr)\operatorname{vec}\bigl(\boldsymbol{\Sigma}^{-1}\bigr).
\end{align*}
\end{itemize}
These and other results are further discussed in~\cite{EntropyBasedStatisticalAnalysisPolSAR}.

Therefore, two additional detectors can be defined by the points that maximize the test statistics based on entropies between the two distributions:
\begin{align*}
\widehat{\jmath}_{\mathcal M}=\arg \displaystyle \max_j 
\underbrace{S_{\mathcal M}(\widehat{\boldsymbol{\theta}_1}{(j)},\widehat{\boldsymbol{\theta}_2}{(N-j)})}_{\equiv \,\,S_{\mathcal M}(j)}
=\arg \displaystyle \max_j S_{\mathcal M}(j),
\end{align*}
where $\mathcal M=\{\text{S},\text{RE-}\beta\}$.

The expressions given in~\eqref{eq:TestStatistic} and~\eqref{eq:EntropyTest} have known asymptotic
distributions so, besides being discriminatory measures, they also
have significance levels. 
This can be useful when the strip of data
conveys a single class and, therefore, no detection should be made.
Nevertheless, in the subsequent section these properties are not
considered, due to the fact that the procedure is bounded to use small
samples.

\section{Applications}~\label{sec:2}

We now apply the methodologies presented in the previous section to simulated and real data. 
Initially, in Section~\ref{application:first}, complex Wishart distributed scenarios were generated in order to quantify and to assess the accuracy and computational load of the detailed edge detection schemes.
For such estimates of the detection error as well as the execution time for each method were compared.
Finally, three applications to real data are performed in Section~\ref{application:second}.
These applications illustrate the difficulties of determining an ideal window size. 

\subsection{Precision and execution time}
\label{application:first}

\begin{figure}[hbt]
\centering
\includegraphics[width=1\linewidth]{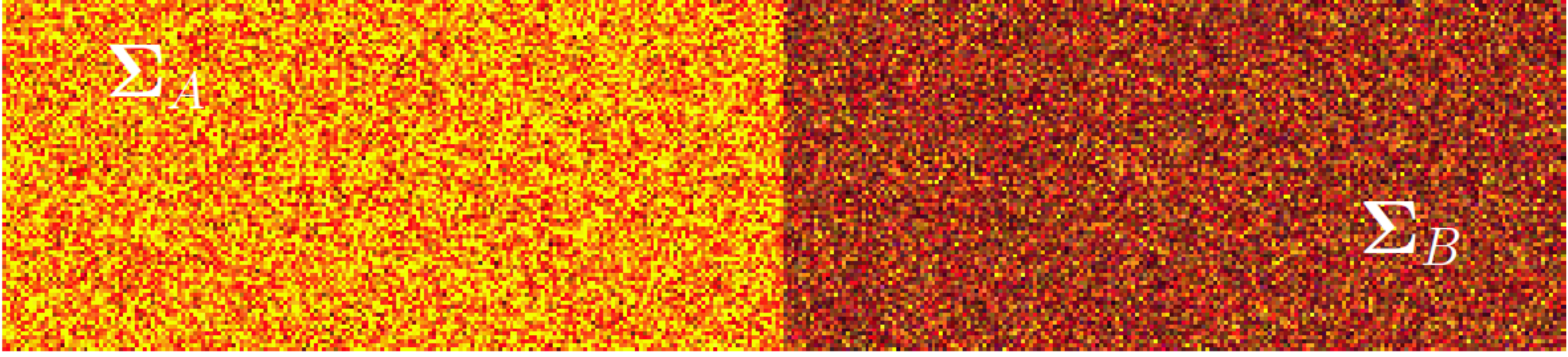}
\caption{
False color of the Pauli decomposition of a simulated image from the scaled Wishart distribution with four looks and two halves defined by $\boldsymbol{\Sigma}_A$ and $\boldsymbol{\Sigma}_B$.
}\label{simuimgf}
\end{figure}

In order to measure and compare the accuracy of the discussed seven edge detection techniques, 
we follow the same methodology proposed by Frery~\textit{et al.}~\cite{PolarimetricSegmentationBSplinesMSSP}. 
We estimate the probability of detecting the edge with an error less than $k$ pixels, 
for $k = 1, 2, \ldots, 10$.
This probability was estimated using strips of data of 400~pixels divided in halves, 
filled with samples from the Wishart distribution with two different covariance matrices:
\begin{align}
\boldsymbol{\Sigma}_A &=\left[\begin{array}{ccc} 
962892 & 19171- 3579\textbf{i} & -154638+191388\textbf{i} \\
& 56707 &   -5798+ 16812\textbf{i}\\
&  &  472251
\end{array} \right] \\
\boldsymbol{\Sigma}_B &=\left[\begin{array}{ccc} 
360932 & 11050+3759\textbf{i} & 63896+1581\textbf{i}\\
& 98960 &   6593+6868\textbf{i}\\
& & 208843
\end{array} \right], 
\label{eq:SigmaB}
\end{align}
which were observed by Frery~\textit{et al.}~\cite{PolarimetricSegmentationBSplinesMSSP} in urban and forest imagery, 
respectively.
In all cases four looks ($L=4$) were considered.
Fig.~\ref{simuimgf} shows a simulated PolSAR image with this configuration according to the Pauli decomposition~\cite{Paulicoding}.

\begin{figure}[hbt]
\centering
\subfigure[Likelihood detector\label{Bias3}]{\includegraphics[width=.49\linewidth]{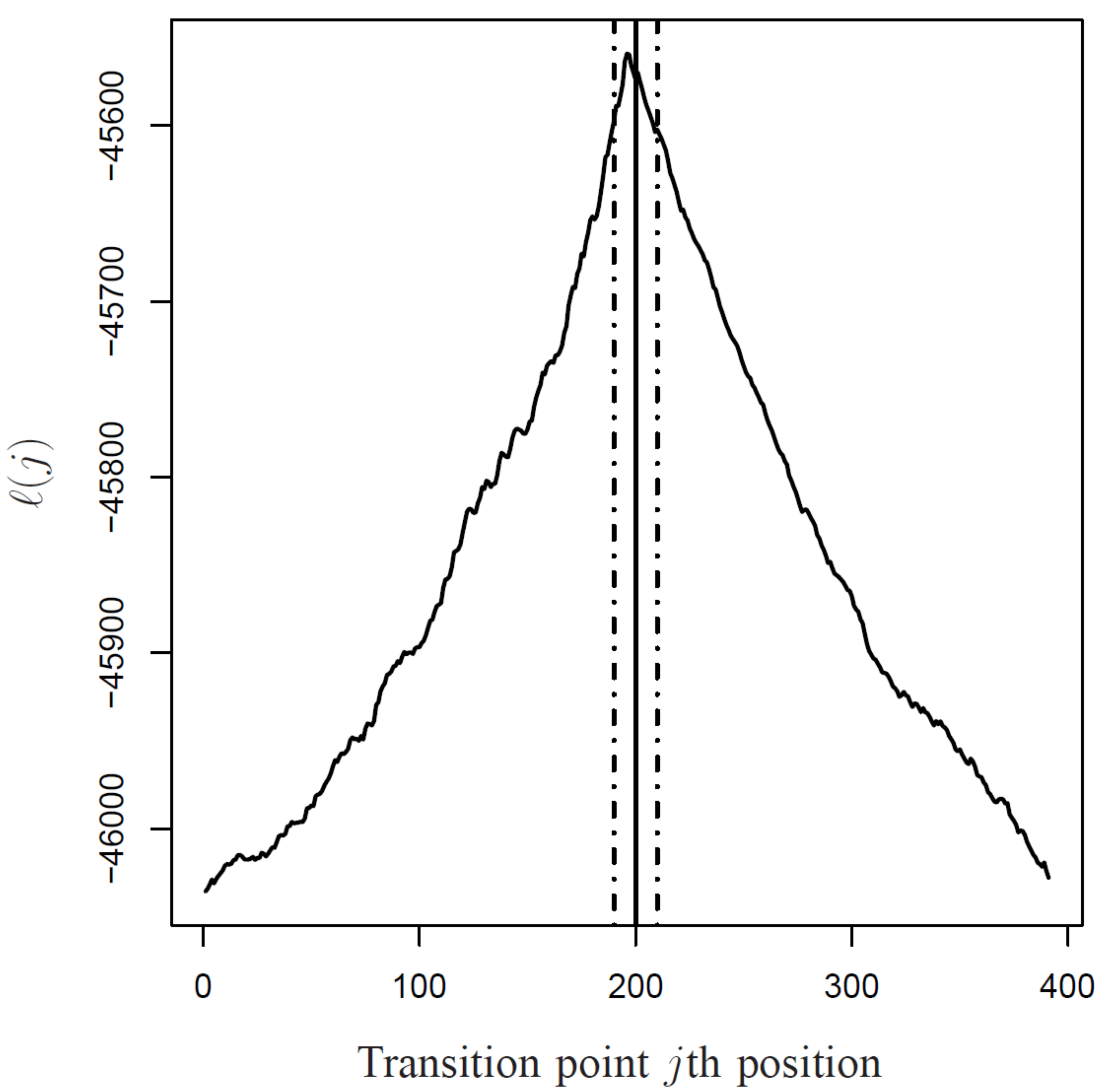}}
\subfigure[Distance detectors\label{Bias1}]{\includegraphics[width=.49\linewidth]{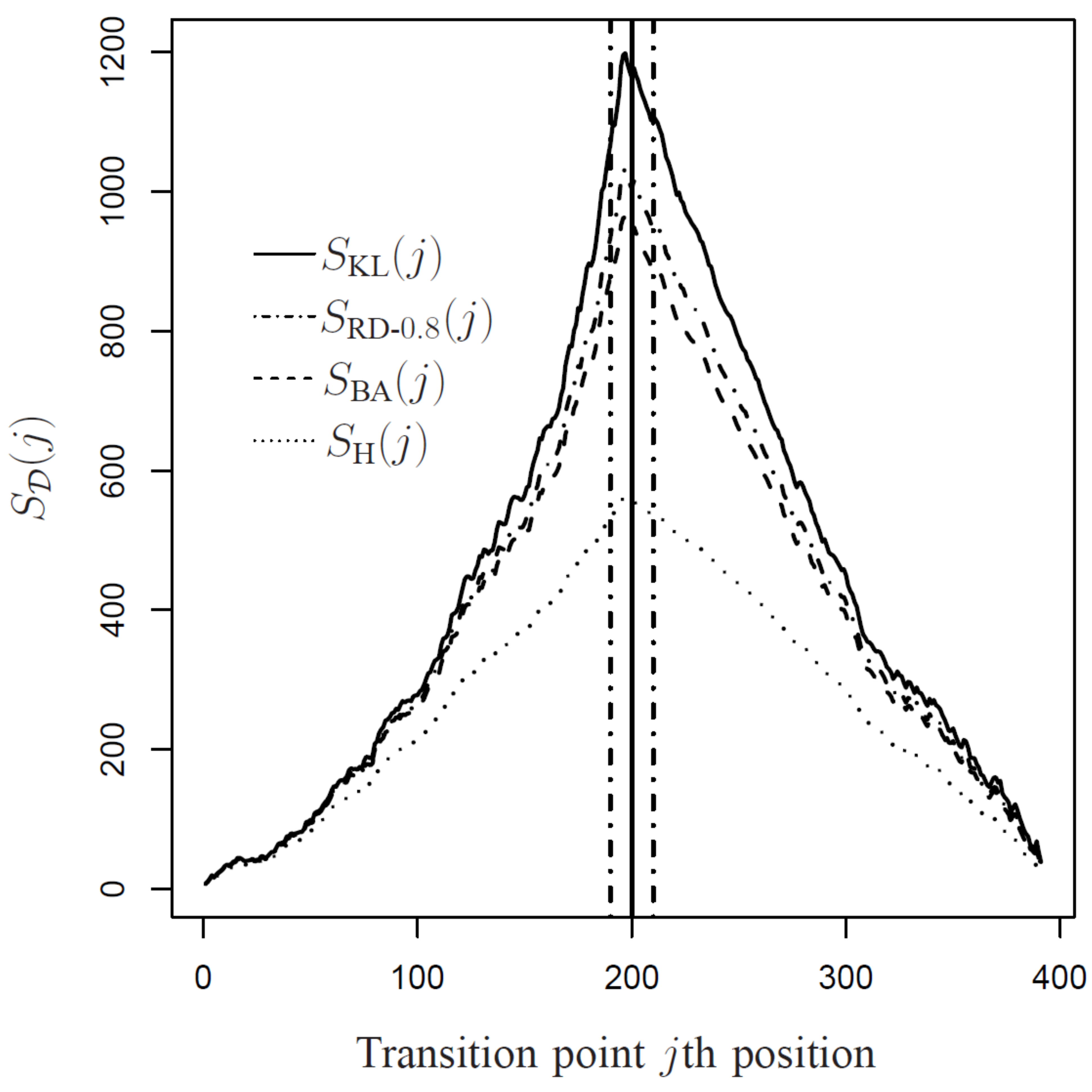}}\\
\subfigure[Entropy detectors\label{Bias22}]{\includegraphics[width=.49\linewidth]{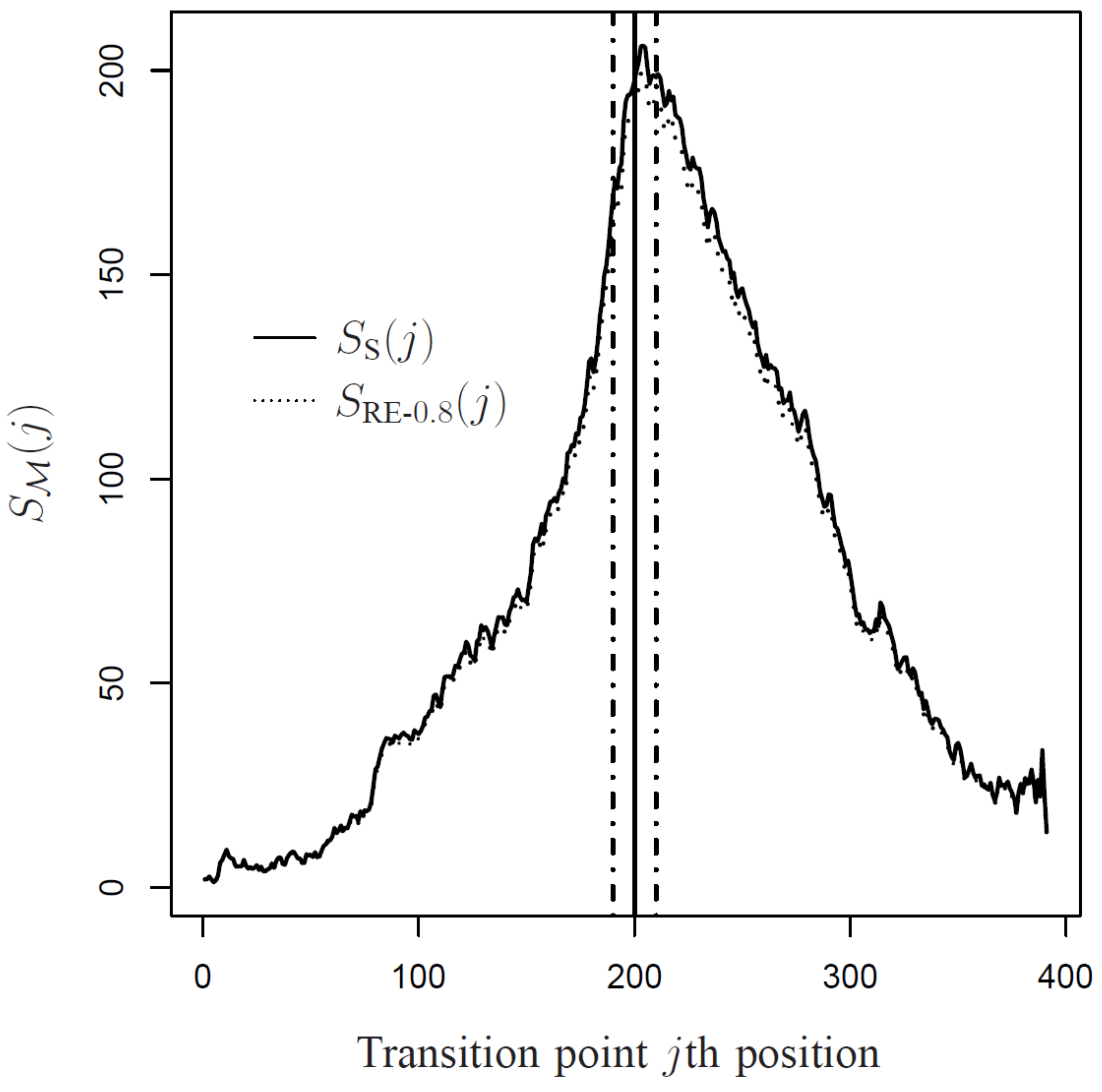}}
\caption{
Illustration of the proposed edge detection measures on images with halves which follow ${\mathcal W}(\boldsymbol{\Sigma_A},4)$ and ${\mathcal W}(\boldsymbol{\Sigma_B},4)$ distributions. The solid vertical line is at the edge and the dashed lines are at $10$ pixels from it.
} 
\label{illustratedmeasures}
\end{figure}

Figure~\ref{illustratedmeasures} shows the functions to be maximized as functions of the position for a typical simulation: the likelihood (Fig.~\ref{Bias3}), distances (Fig.~\ref{Bias1}), and entropies (Fig.~\ref{Bias22}).
The edge is at 200, and this position identified with a vertical solid line; the symmetric interval of twenty pixels around it is shown in vertical dash-dot lines.
All of these distances exhibit a maximum which is close to the true edge position.

We generated 1000~independent scene simulations to obtain the estimated boundaries positions 
($b(r)\in\{\widehat{\jmath}_{\text{ML}}$, $\widehat{\jmath}_{\mathcal D}$, $\widehat{\jmath}_{\mathcal M}\}$, $1\leq r\leq 1000$).
The distance between these points and the true boundary (the absolute empirical bias) was evaluated for each replication:
$$
E(r)=|b(r)-200|,
$$
and the probability of observing an error smaller than a certain number of pixels is estimated by relative frequencies as 
$$
f(k)=T(k)/1000,
$$
where $T(k)$ the number of replications for which the error is smaller than $k$ pixels.

\begin{figure}[hbt]
\centering
\includegraphics[width=1\linewidth]{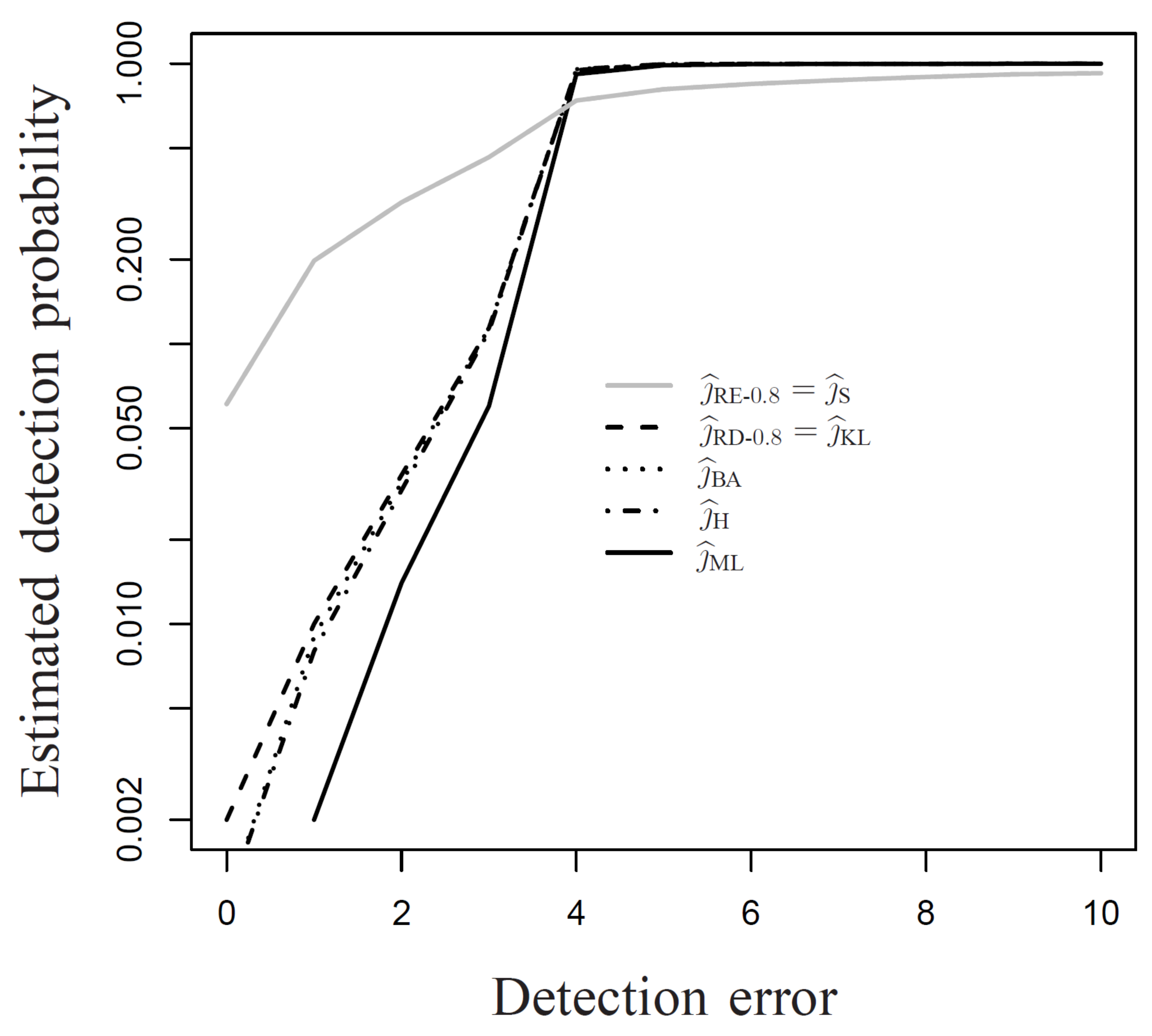}
\caption{
Performance of detectors on phantom images with two halves distributed according to scaled complex Wishart laws.
} 
\label{Figure1}
\end{figure}

Fig.~\ref{Figure1} presents, in semilogarithmic scale, the estimated probability of finding the edge with an error equal or smaller than $k$ pixels, $k=1,2,\ldots,10$.
Results show evidence that the methods based on entropies have the best results for $k\leq 3$.
On the other hand, admitting an error $k\geq 4$, the method involving stochastic distances have the higher estimated detection probability.
This fact suggests that a combination of these methodologies could lead to even more precise edge detection.
Finally, the method based on the likelihood is outperformed by those based on stochastic distances; however, the former method provides  better estimates when compared to $S_{\text{S}}$ and $S_{\text{R-}\beta}$.

\subsubsection{Influence of the backscatter matrix and time execution}

\begin{figure}[hbt]
\centering
\includegraphics[width=1\linewidth]{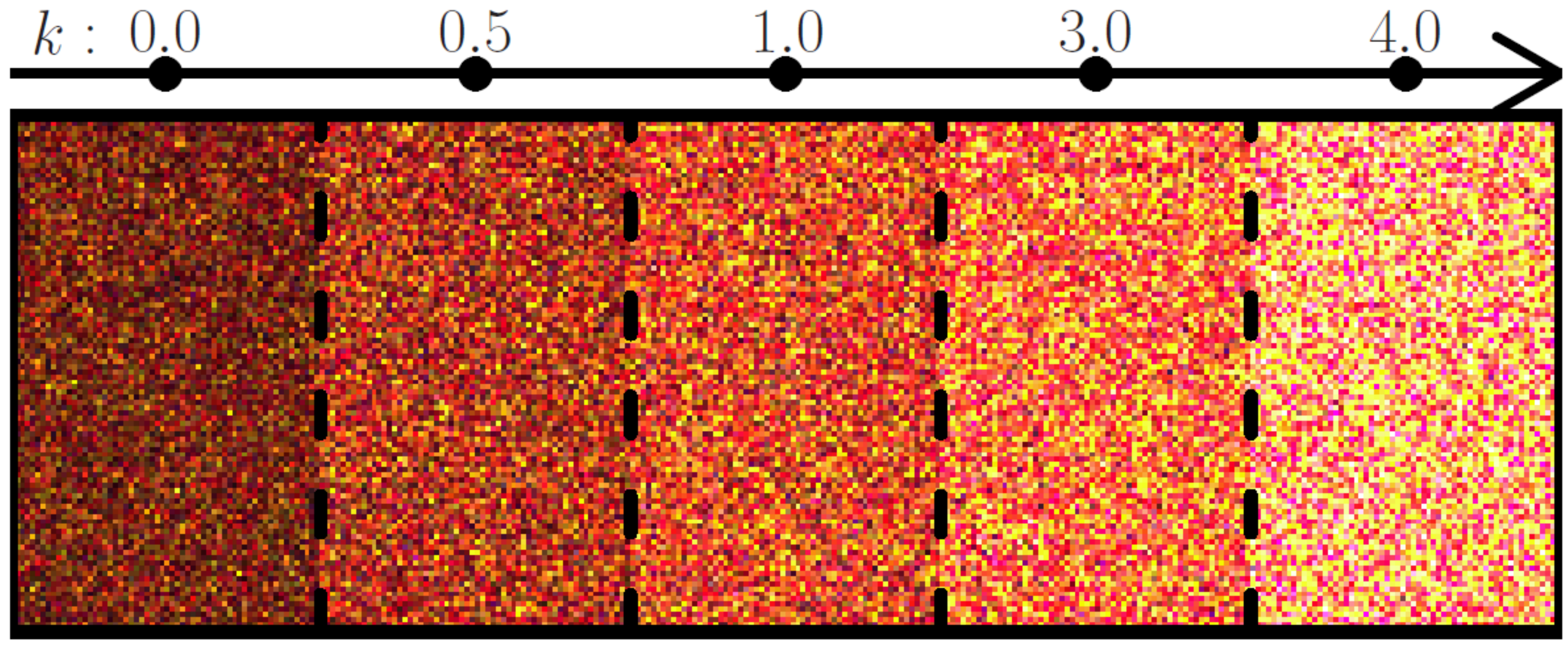}
\caption{
Images generated from $W(\boldsymbol{\Sigma_B} \cdot (1 + k); 4)$, for $k = 0, 1/2, 1, 3, 4$, mapping [HH$+$VV] -- [HV] -- [HH$-$VV] onto the RGB channels.
}
\label{Figure2}
\end{figure}

Now we devise a simulation study
where data is more similarly distributed;
therefore offering a more difficult scenario.
Also we aim at assessing
the execution time of each method.
Strips of data were defined with halves 
which follow distributions ${\mathcal W}((1+k)\boldsymbol{\Sigma}_B,4)$ and ${\mathcal W}(\boldsymbol{\Sigma}_B,(1+v)4)$,
respectively. 
This situation is referenced by ``case-$(k,v)$''.
Fig.~\ref{Figure2} presents images sampled in cases $(0,0)$, $(1/2,0)$, $(1,0)$, $(3,0)$, and $(4,0)$, i.e., with the same number of looks and varying covariance matrices.

\begin{figure}[hbt]
\centering
\includegraphics[width=1\linewidth]{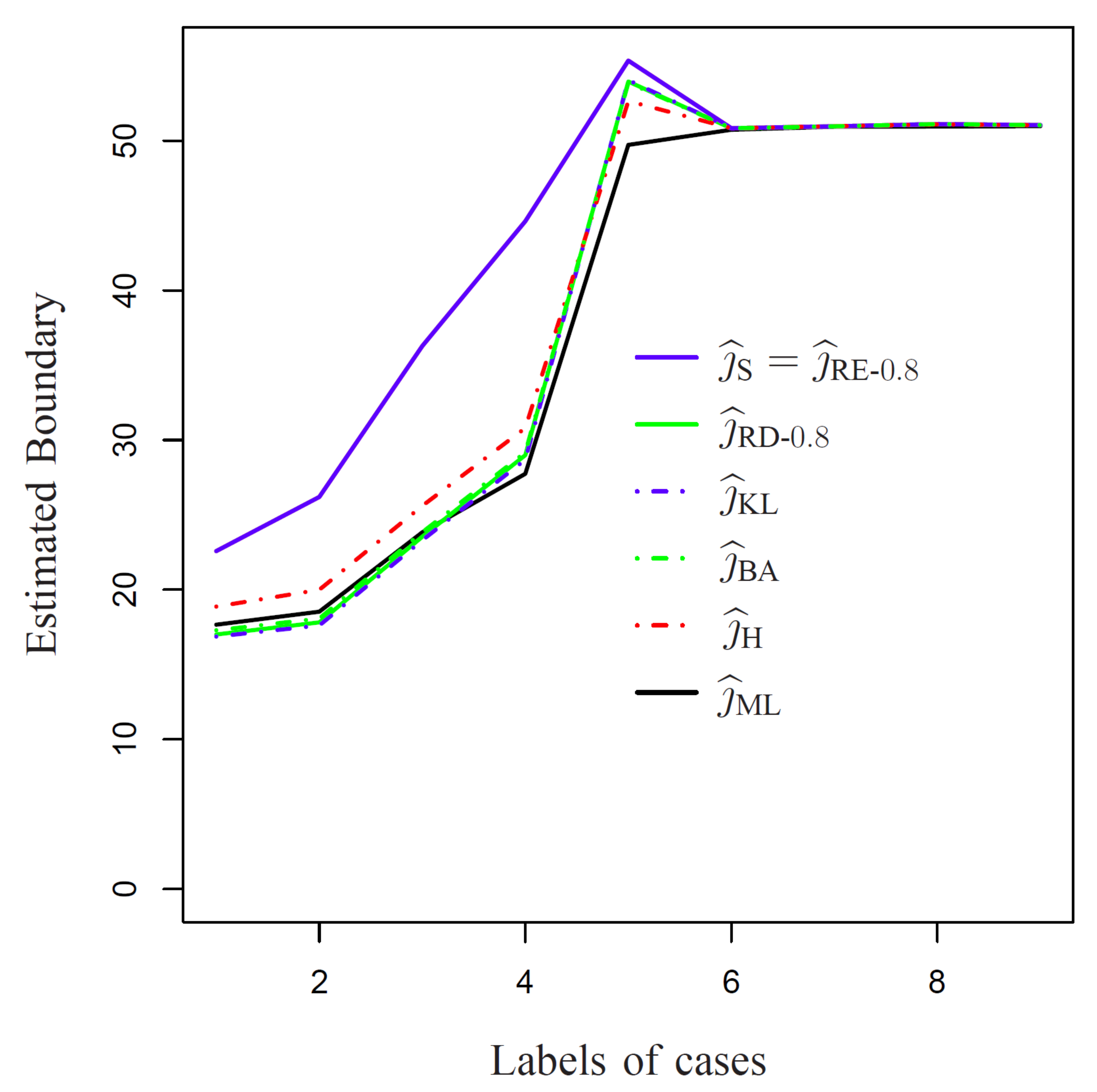}
\caption{Precison of the detection on polarimetric images with halves which follow ${\mathcal W}((1+k)\boldsymbol{\Sigma}_B,4)$ and ${\mathcal W}(\boldsymbol{\Sigma}_B,(1+v)4)$ distributions.
} 
\label{simulationstudyII}
\end{figure}

The precision and execution time of the procedure 
were estimated using
one thousand replications. 
The latter 
was measured 
in seconds of CPU cycles, while for the former we adopted a ``hit-and-miss'' criterion:
a run is considered successful if it finds the true boundary,
otherwise it is assumed wrong. 
The mean number of successful runs 
are reported 
in Fig.~\ref{simulationstudyII}.

\begin{figure}[hbt]
\centering
\subfigure[All times\label{time3}]{\includegraphics[width=.48\linewidth]{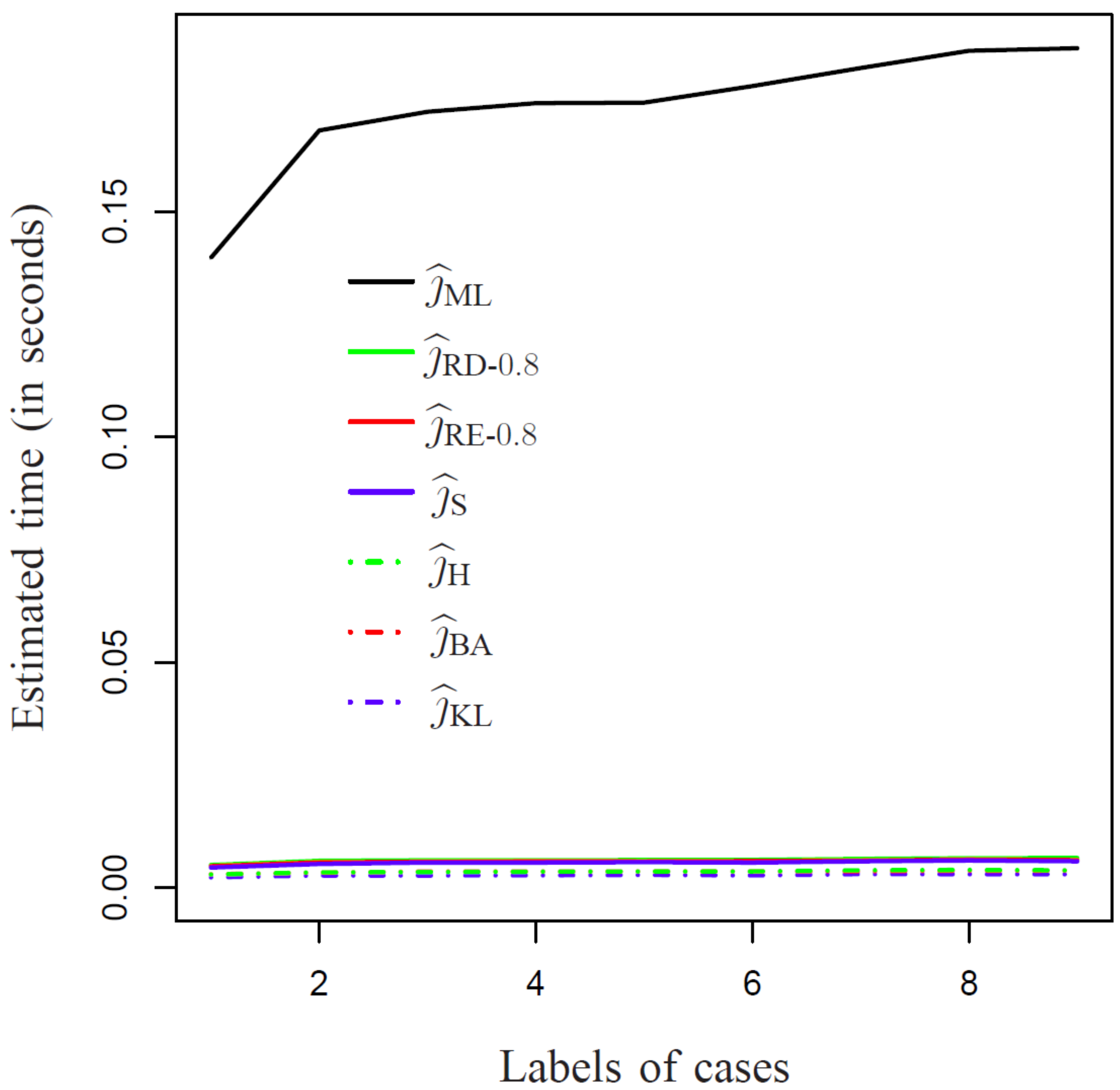}}
\subfigure[Information theoretic tool times\label{time1}]{\includegraphics[width=.48\linewidth]{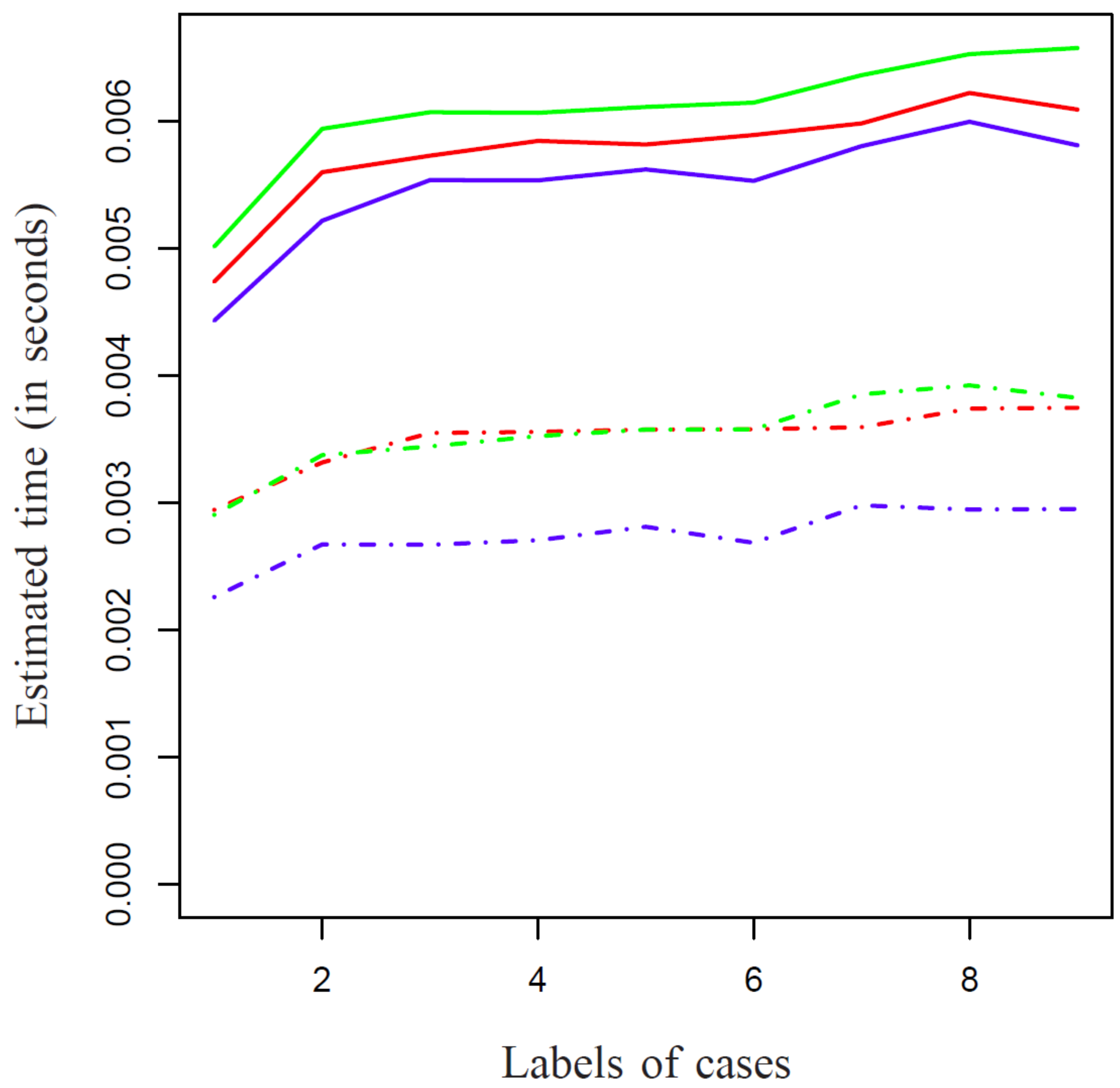}}
\caption{Detection times of the proposed methods in terms of the situations presented in Table~\ref{Estimate1}.
}
\label{timestudy}
\end{figure}

Fig.~\ref{simulationstudyII} and Table~\ref{Estimate1} present the average estimates of precision for all methods 
under several situations; the first column of Table~\ref{Estimate1} shows the label for each case which is used in Figs.~\ref{simulationstudyII} and~\ref{timestudy}.
The best estimators are highlighted in boldface.

Situations $1$ to $4$, which correspond to cases where the covariance matrix does not vary or changes by a small amount, 
the hardest ones to perform edge detection.
Thus they often lead to low precision regardless the technique.
It is noticeable that entropy methods consistently outperform other techniques in these particularly challenging situations.
The performance of all procedures is comparable at remaining situations.
We could identify that the maximum likelihood and entropy-based  estimators
excelled, in three and in two situations, respectively, but for a small difference of the value of the true boundary.

\begin{table}[hbt]                                                     
\centering       
\setlength{\tabcolsep}{2pt}                                                                                                                  
\caption{Estimates for boundary points (estimates closest to the true edge are highlighted)}\label{Estimate1}                
\begin{tabular}{*8{c}}\toprule  
Label & Case  & $\widehat{\jmath}_\text{ML}$ & $\widehat{\jmath}_\text{KL}$ & $\widehat{\jmath}_\text{BA}$ & $\widehat{\jmath}_\text{H}$ &    
 $\widehat{\jmath}_{\text{RD-}0.8}$ & $\widehat{\jmath}_\text{S}=\widehat{\jmath}_{\text{RE-}0.8}$ \\                                                                                             
\cmidrule(lr{.25em}){1-2}\cmidrule(lr{.25em}){3-8}                      
$1$ & $(0,3)$  & $17.662$ & $16.870$ & $17.288$ & $18.871$ & $17.007$ & $\boldsymbol{22.582}$  \\
$2$ & $(0,7)$  & $18.531$ & $17.609$ & $18.074$ & $19.994$ & $17.816$ & $\boldsymbol{26.198}$  \\
$3$ & $(0,1)$  & $23.834$ & $23.306$ & $23.843$ & $25.604$ & $23.538$ & $\boldsymbol{36.287}$  \\
$4$ & $(0.1,1)$& $27.746$ & $28.688$ & $29.176$ & $30.806$ & $28.991$ & $\boldsymbol{44.626}$  \\
$5$ & $(0.1,0)$& $\boldsymbol{49.726}$ & $54.055$ & $53.877$ & $52.653$ & $53.961$ &  $55.358$  \\
$6$ & $(1,3)$  & $50.738$ & $50.826$ & $50.841$ & $50.842$ & $50.834$ & $\boldsymbol{50.844}$  \\
$7$ & $(2,7)$  & $50.941$ & $50.966$ & $50.971$ & $50.974$ & $50.968$ & $\boldsymbol{50.970}$  \\
$8$ & $(1,0)$  & $\boldsymbol{50.952}$ & $51.129$ & $51.119$ & $51.102$ & $51.120$ & $51.110$  \\
$9$ & $(2,0)$  & $\boldsymbol{50.979}$ & $51.052$ & $51.048$ & $51.042$ & $51.050$ & $51.052$  \\                                                  
\bottomrule                                                             
\end{tabular}                                                           
\end{table}                                                            

Fig.~\ref{timestudy} presents the mean execution time of each method in seconds.
Fig.~\ref{time3} presents obtained detection times, where it is noticeable that the
maximum likelihood method associated to $\widehat{\jmath}_\text{ML}$
is the slowest procedure.
Fig.~\ref{time1} shows only the times required by procedures based on stochastic distances.
We notice that the 
method based on the Kullback-Leibler divergence,
which furnishes $\widehat{\jmath}_\text{KL}$,
was the fastest one.
The results can be summarized in the following inequalities of the required computational times $t$ for each technique:
$$
t_{\widehat{\jmath}_\text{ML}} \geq  t_{\widehat{\jmath}_\text{RD}} \geq t_{\widehat{\jmath}_\text{RE}} \geq t_{\widehat{\jmath}_\text{S}} \geq t_{\widehat{\jmath}_\text{H}} \geq t_{\widehat{\jmath}_\text{BA}} \geq  t_{\widehat{\jmath}_\text{KL}}.
$$ 
The influence of the case on the computing time is very small.
All implementations are in the \texttt{R} version 2.13.2 programming language, running on a PC with an 
{Intel\textcopyright\ Core\textregistered\ i7-2630QM} \unit{2.00}{GHz} with \unit[4]{GB} RAM.

\subsubsection{Influence of the backscatter matrix SPAN and image resolution}

In the following we assess the relative information provided by the polarimetric data with respect to single-channel data for the problem of finding edges.
Assuming the scaled complex Wishart with parameters $\boldsymbol{\Sigma}$ and $L$ for the polarimetric data, each intensity channel is described by the gamma distribution with density given by
$$
f_{Z_i}(Z'_i;L/\sigma^2_{i},L)=\frac{L^L{Z'_i}^{L-1}}{ \sigma^{2L}_i\Gamma(L)}
\exp\bigl(-LZ'_i/\sigma^2_i\bigr),
$$
for $i\in\{\text{HH,HV,VV}\}$, 
where $\sigma^2_i$ is the $(i,i)$th entry of $\boldsymbol{\Sigma}$
and $Z'_i$ is the $(i,i)$th entry of the random matrix $\boldsymbol{Z}$~\cite{Hagedorn2006655}.  
This density can be used in the log-likelihood of strips (cf.~\eqref{eq:loglik}) to derive the estimator of the edge position (cf.~\eqref{eq:jML}).

We considered two random matrices $\boldsymbol{X}\sim{\mathcal W}(\boldsymbol{\Sigma}_B,4)$ and $\boldsymbol{Y}\sim{\mathcal W}(\boldsymbol{\Sigma}_{B'},4)$:
\begin{itemize}
\item Case A: $\operatorname{diag}(\boldsymbol{\Sigma}_{B'})=(1+\delta) \operatorname{diag}(\boldsymbol{\Sigma}_B)$ and $\operatorname{adiag}(\boldsymbol{\Sigma}_{B'})=\operatorname{adiag}(\boldsymbol{\Sigma}_B)$, 
where $\operatorname{diag}(\boldsymbol{A})=[A_{11}\,A_{22}\,\ldots\,A_{mm}]^\top$, for a square matrix $\boldsymbol{A}$ with order $m$, and $\operatorname{adiag}(\boldsymbol{A})=[\boldsymbol{w}^\top\,\boldsymbol{w}^*]^\top$, $\boldsymbol{w}=[A_{12}\,$ $\ldots\,A_{1m}\,$ $A_{23}\,$ $\ldots\,$ $A_{2m}\,$ $\ldots\,$ $A_{(m-1)m}]^\top$.
This case represents the influence of the SPAN~\cite{CaoHongWuPottier2007}, i.e., the trace of backscatter matrix when $\delta\geq 0$.
\item Case B: $\operatorname{diag}(\boldsymbol{\Sigma}_{B'})=\operatorname{diag}(\boldsymbol{\Sigma}_B)$ 
and
$\operatorname{adiag}(\boldsymbol{\Sigma}_{B'})=\operatorname{adiag}(\boldsymbol{\Sigma}_B)+(1+\delta) \{\Re[\operatorname{adiag}(\boldsymbol{\Sigma}_B)]+\Im[\operatorname{adiag}(\boldsymbol{\Sigma}_B)]\}$. 
This aims at quantifying the influence of the covariance between polarization channels varying the values of $\delta\geq -1$.
\end{itemize}
A Monte Carlo simulation was performed to quantify the performance of the edge detectors under these conditions.

\begin{figure}[hbt]
\centering
\subfigure[$\delta=0.2$\label{parsim1}]{\includegraphics[width=.48\linewidth]{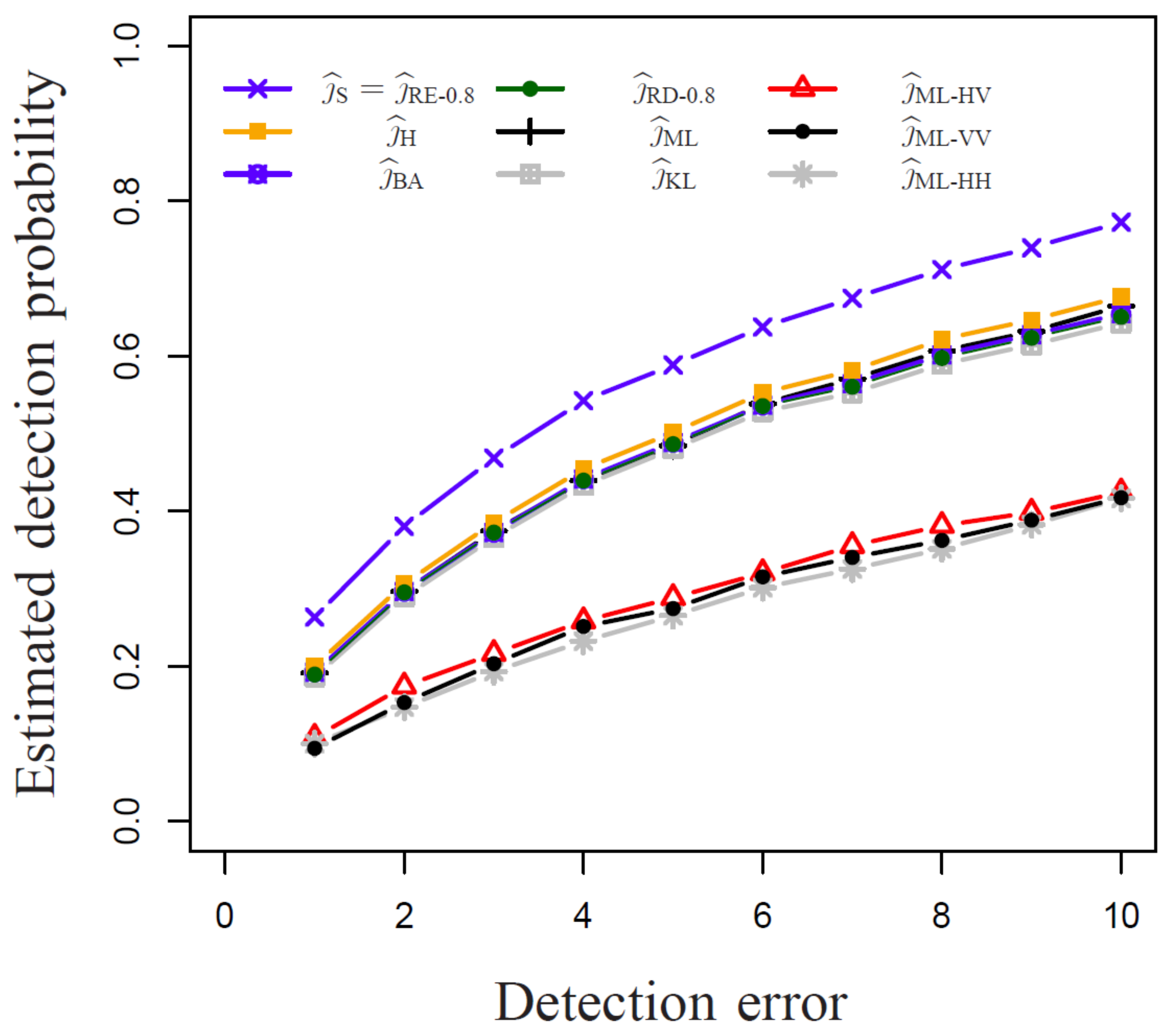}}
\subfigure[$\delta=0.4$\label{parsim2}]{\includegraphics[width=.48\linewidth]{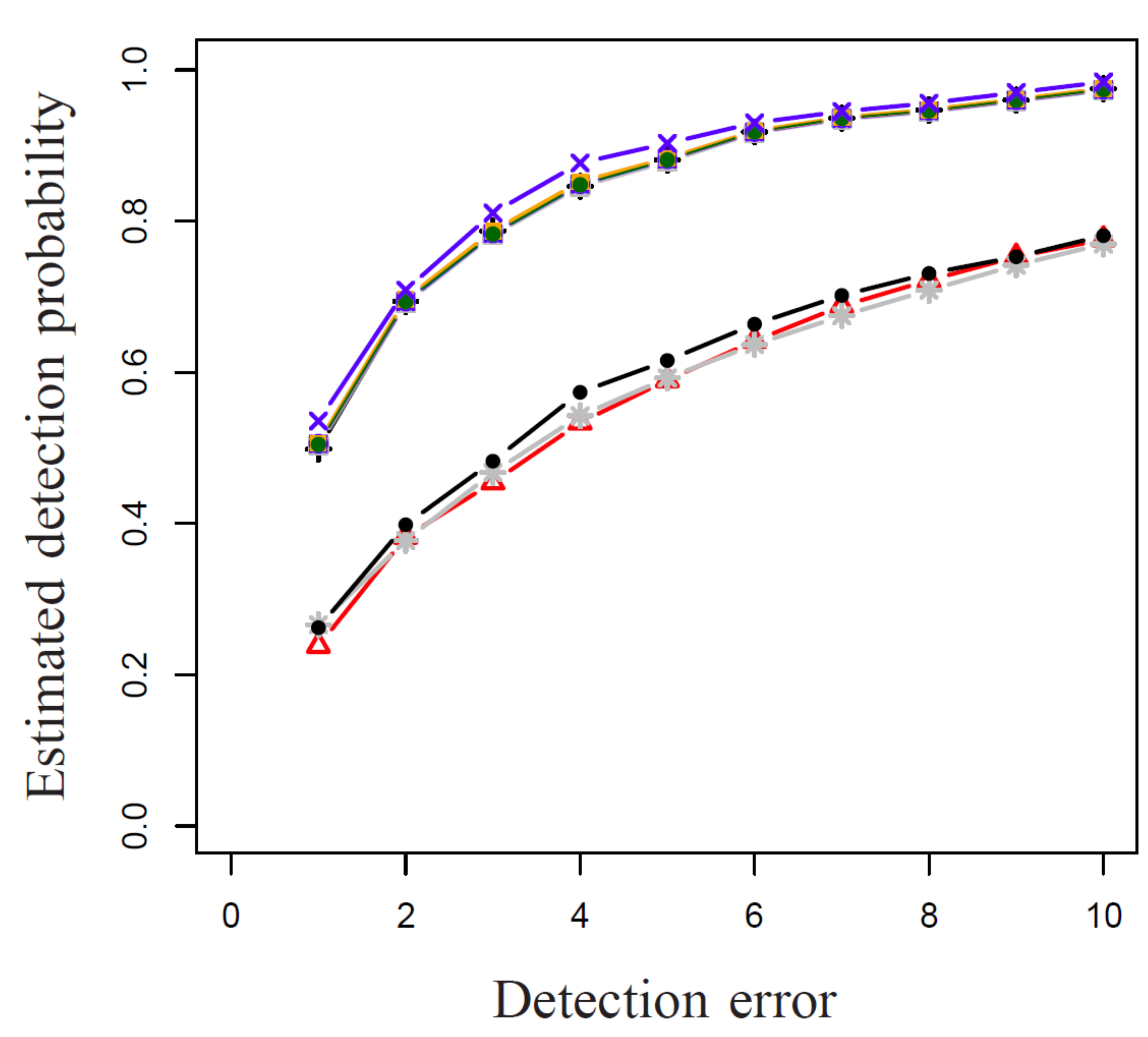}}\\
\subfigure[$\delta=1.0$\label{parsim3}]{\includegraphics[width=.48\linewidth]{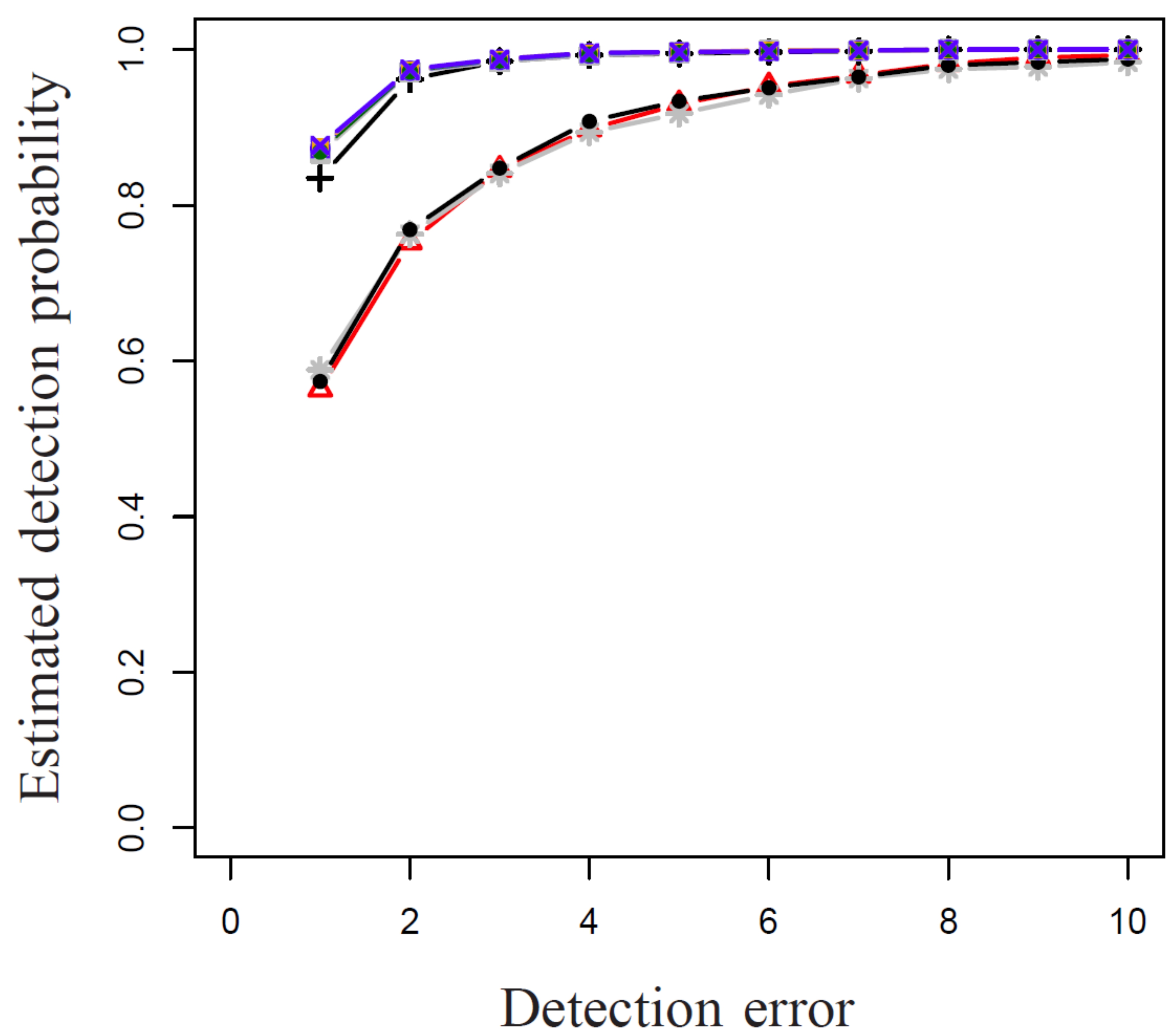}}
\caption{
Influence of the SPAN on the probability of detecting and edge with an error equal or smaller than a certain number of pixels.
} 
\label{Simulationdiagonal}
\end{figure}

For Case~A, $1000$ images with $100$ pixels in each half were generated with $\delta\in\{0.2,0.4,1.0\}$.
Edge detection was performed, and the probability of detecting within an error equal or smaller than 
$k\in\{1,2,\ldots,10\}$ pixels was estimated.
Fig.~\ref{Simulationdiagonal} presents the results.
As expected, all polarimetric methods outperformed 
the methods based on single-channel data.
Additionally, 
these results provide evidence that the performance of the methods tends to be similar when the SPAN is significantly different, i.e., when the halves are quite different.
Moreover,
boundary detection based on entropy measures outperformed 
all other methods in Case~A.

\begin{figure}[hbt]
\centering
\subfigure[$\delta=0.2$\label{parsim11}]{\includegraphics[width=.48\linewidth]{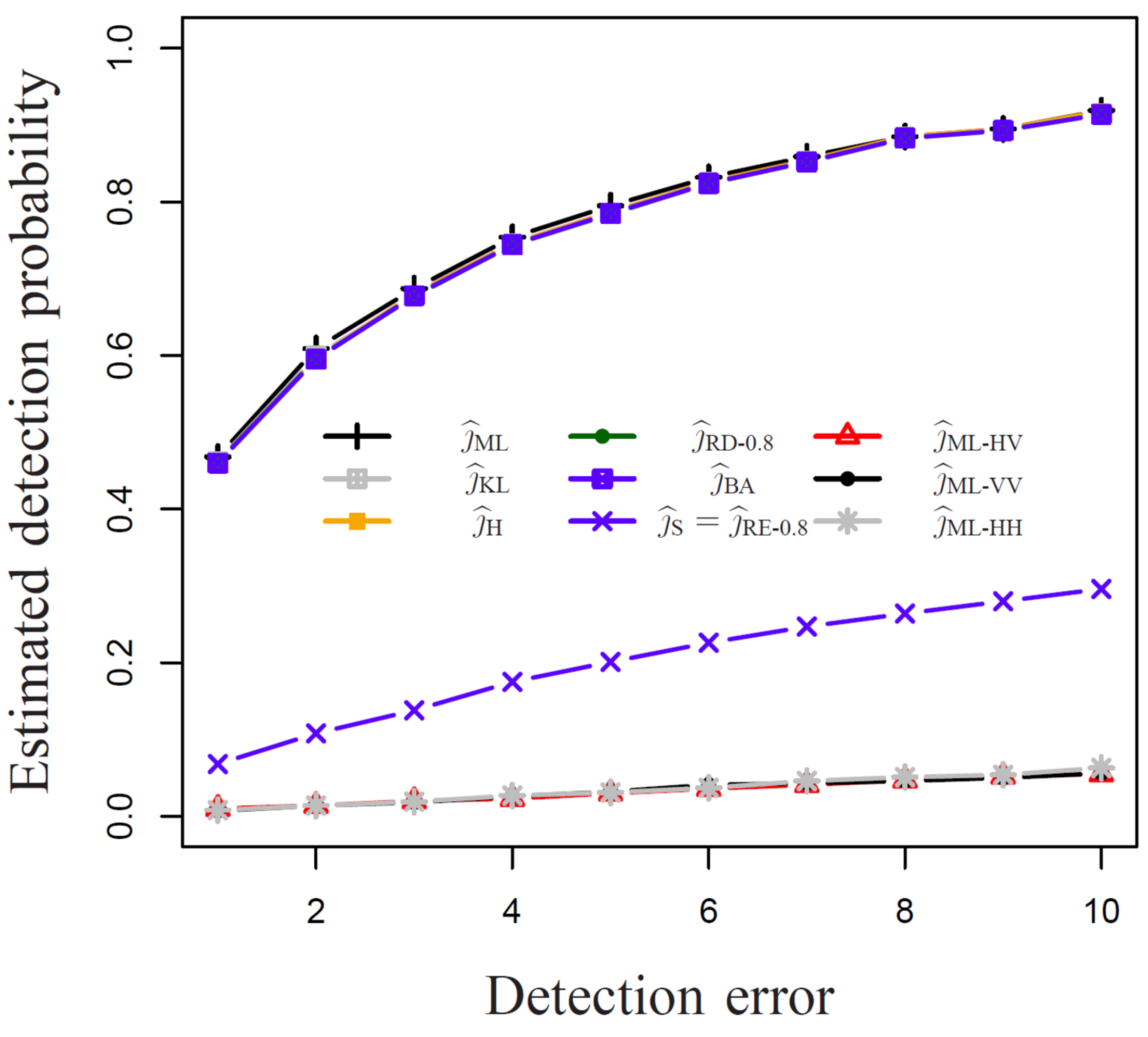}}
\subfigure[$\delta=0.4$\label{parsim22}]{\includegraphics[width=.48\linewidth]{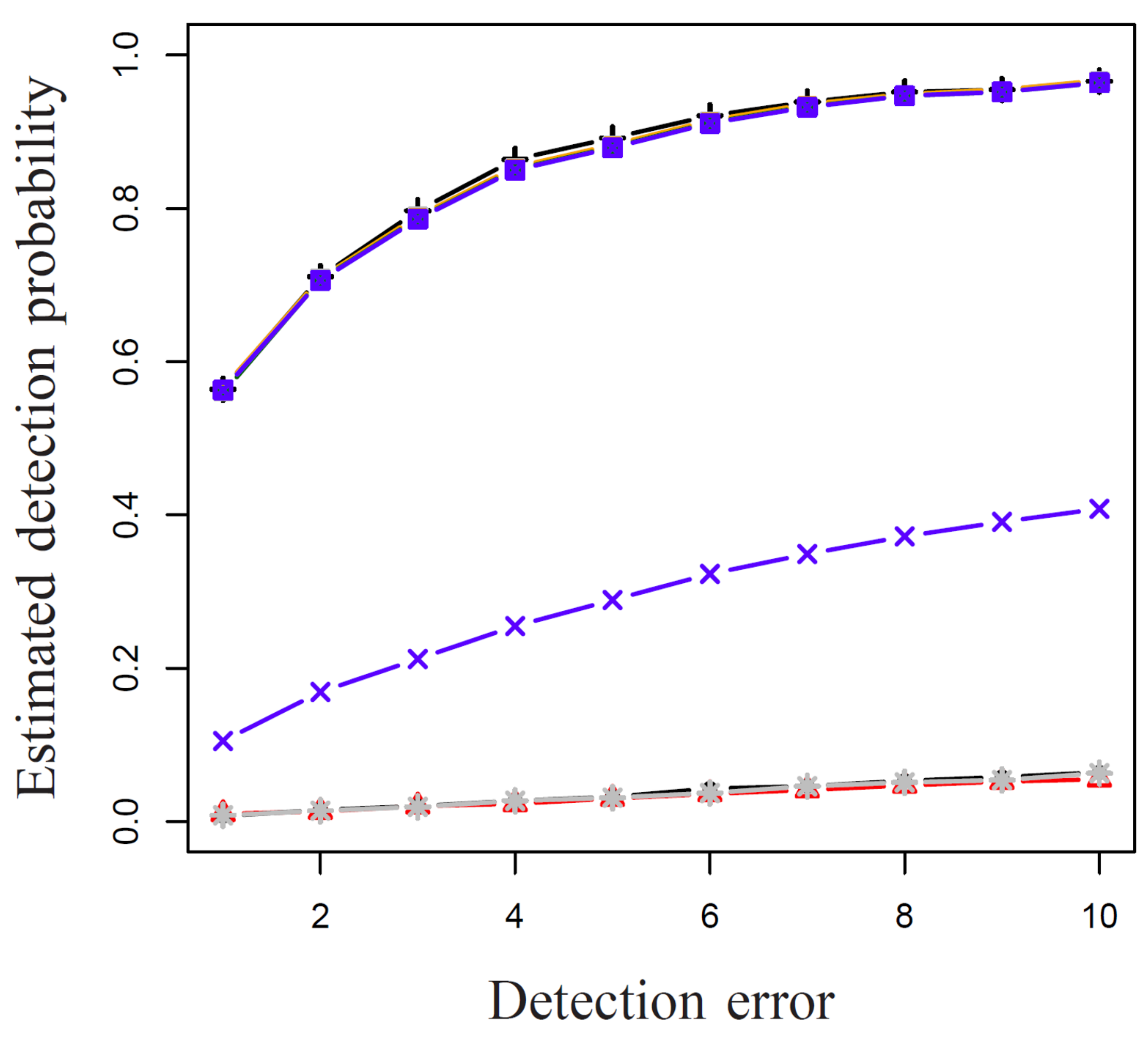}}\\
\subfigure[$\delta=1.0$\label{parsim33}]{\includegraphics[width=.48\linewidth]{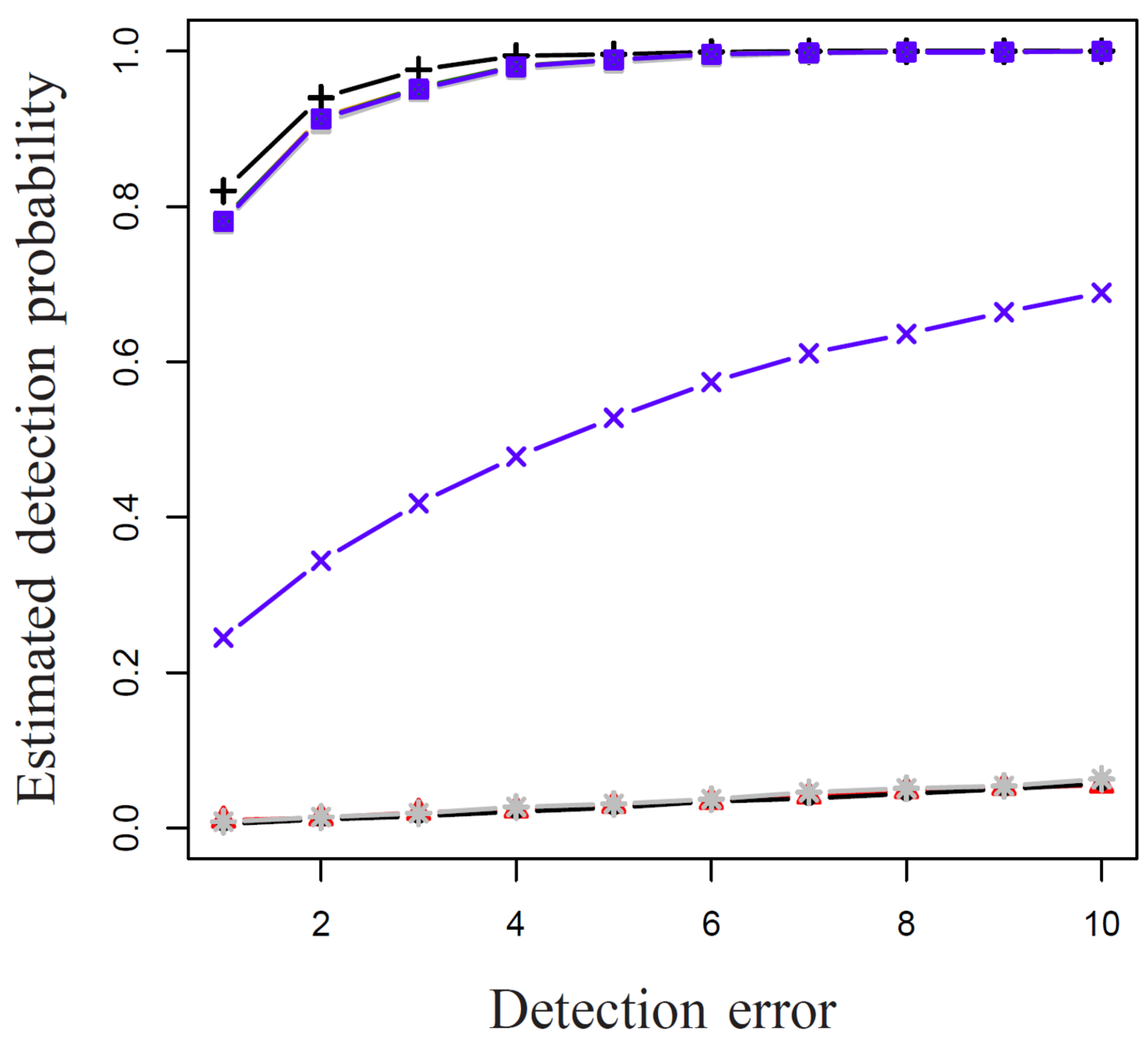}}
\caption{
Influence of the covariance between polarization channels on the probability of detecting and edge with an error equal or smaller than a certain number of pixels.
} 
\label{SimulationdiagonalI}
\end{figure}

Fig.~\ref{SimulationdiagonalI} shows the performance of the discussed methods in Case~B, for $\delta\in\{0.2,0.3,1.0\}$.
As in Case~A, polarimetric data led to better results than single-channel images.
In this case, detectors based on stochastic distances and maximum likelihood outperformed those based on entropies. 
Fig.~\ref{parsim11} legend presents the detectors ordered from best to worst from top to bottom and from left to right.
The best result was obtained by $\widehat{\jmath}_{\text{ML}}$, but small differences were observed among methods involving distances. 

We now analyze the effect of spatial resolution on the precision of the edge detection procedures.
To that end, we simulated a strip of 200~observations divided in two halves, each with different distributions, namely $\boldsymbol{X}\sim{\mathcal W}(\boldsymbol{\Sigma}_B,4)$ and $\boldsymbol{Y}\sim{\mathcal W}(\boldsymbol{\Sigma}_{B'},4)$, such that $\operatorname{diag}(\boldsymbol{\Sigma}_{B'}) = 1.2 \operatorname{diag}(\boldsymbol{\Sigma}_B)$, with $\boldsymbol{\Sigma}_B$ defined in~\eqref{eq:SigmaB}.
This strip is assumed to have the best spatial resolution, and is denoted ``$1\div 1$''.
The position of the edge, which is at the middle of the strip, i.e. at 100, is estimated with all the available techniques.
Then the spatial resolution of the strip is degraded by taking the mean of pairs of contiguous observations and downsampling one each two pixels obtaining, thus, a $1\div2$ resolution image.
The estimation of the edge position, which is now at 50, is performed.
A new pyramid is formed, denoted ``$1\div4$'', and the estimation of the edge (now at $25$) is again performed.
The equivalent number of looks $L$ is assumed know, and at each resolution takes the values $4$, $8$ and $16$.
This procedure was repeated generating $1000$ independent initial strips, and the bias ($\widehat B$, the expected value of the estimator minus the true value), the standard deviation ($\operatorname{sd}$), the coefficient of variation ($\operatorname{CV}$) and the mean squared error ($\operatorname{MSE}$) of each estimator $\widehat{\jmath}_\bullet)$ were estimated.
The results are shown in Table~\ref{Estimate2}.

The detection bias was always very small, being the largest values
in the order of $10^{-2}$ of a pixel and 
the smallest ones in the order of $10^{-3}$ of a pixel.
The absolute bias tended to exhibit a small increase when the spatial resolution was reduced, while the standard deviation and the mean squared error were reduced.
The width of the Gaussian symmetric confidence intervals at the $95\%$ confidence level ranges between, approximately, half a pixel and one pixel, rendering, thus, that all techniques exhibit comparable accuracy.
Changing the resolution had little effect on the coefficient of variation.

\begin{table*}[hbt]
\setlength{\tabcolsep}{3pt}
\centering
\caption{Effect of the spatial resolution on the edge detection}\label{Estimate2}
\begin{tabular}{*2{c} *9{r}}
\toprule
& \multicolumn{9}{c}{Boundary detection methods} \\ \cmidrule(lr{.25em}){3-11}
& & \multicolumn{3}{c}{Single channel data} & \multicolumn{6}{c}{Full PolSAR data} \\ \cmidrule(lr{.25em}){3-5}\cmidrule(lr{.25em}){6-11}
Resolution  & Measures &  $\widehat{\jmath}_\text{HH-ML}$ & $\widehat{\jmath}_\text{HV-ML}$ & $\widehat{\jmath}_\text{VV-ML}$ & 
$\widehat{\jmath}_\text{ML}$ & $\widehat{\jmath}_\text{KL}$ & $\widehat{\jmath}_\text{BA}$ & $\widehat{\jmath}_\text{H}$ & 
$\widehat{\jmath}_{\text{RD-}0.8}$ & $\widehat{\jmath}_\text{S}=\widehat{\jmath}_{\text{RE-}0.8}$ \\ 
\cmidrule(lr{.25em}){1-2}\cmidrule(lr{.25em}){3-11}
$1\div 1$ & $100 \widehat{B}$ & $0.502$ & $-0.547$ & $2.866$ & $0.273$ & $-0.793$ & $-0.785$ & $-0.545$ & $-0.792$ & $-0.171$ \\
& $\operatorname{sd}(\widehat{\jmath}_{\bullet})$ & 
 $ 52.850$ &   $48.948$ &  $50.683$ & $18.388$ & $24.338$ & $22.733$ & $18.826$ & $24.338$ & $15.028$  \\
& $\operatorname{CV}(\widehat{\jmath}_{\bullet})$ & 
 $  0.263$ &    $0.246$ &    $0.246$ &   $0.092$ &   $0.123$ &  $0.115$ & $0.095$ & $0.123$ & $0.075$  \\
& $\operatorname{MSE}(\widehat{\jmath}_{\bullet})$ &
 $2791.361$ & $2394.661$ & $2599.036$ & $338.076$ & $594.280$ & $518.758$ & $355.249$ & $594.232$ & $225.726$ \\
 \cmidrule(lr{.25em}){2-11}
$1\div 2$ & &  $1.541$ & $-0.252$ & $3.335$ & $0.808$ & $-0.151$ & $-0.153$ & $-0.052$ & $-0.152$ & $0.420$ \\
&&            $ 25.435$  & $ 22.977$  & $ 25.061$  & $  8.984$ & $  9.880$ & $  9.875$ & $  9.406$ & $  9.875$ & $  7.373$  \\
&&            $  0.250$  & $  0.230$  & $  0.243$  & $  0.089$ & $  0.099$ & $  0.099$ & $  0.094$ & $  0.099$ & $  0.073$  \\
&&            $648.653$  & $527.482$  & $638.523$  & $ 81.288$ & $ 97.549$ & $ 97.433$ & $ 88.386$ & $ 97.434$ & $ 54.486$  \\
\cmidrule(lr{.25em}){2-11}
$1\div 4$ & & $2.144$ & $1.190$ & $4.504$ & $1.876$ & $0.784$ & $0.932$ & $0.990$ & $0.890$ & $1.362$ \\
&&            $ 12.037$  & $ 11.235$ & $ 11.790$ & $ 4.451$ & $ 4.933$ & $ 4.713$ & $ 4.671$ & $ 4.737$ & $ 3.603$ \\
&&            $  0.236$  & $  0.222$ & $  0.226$ & $ 0.087$ & $ 0.098$ & $ 0.093$ & $ 0.092$ & $ 0.094$ & $ 0.071$ \\
&&            $145.882$  & $126.445$ & $143.948$ & $20.670$ & $24.468$ & $22.404$ & $22.039$ & $22.617$ & $13.435$ \\
\bottomrule                                                                                                                
\end{tabular}                                                                                                              
\end{table*}                                                                                                                

\subsection{Application to real data}\label{application:second}

We applied the methods to three real images. 
The first example describes a simple situation: well separated regions with most of the strips including two regions only. 
The second and third examples present more sophisticated problems due to the size of the strips (larger strips may span more than two regions).

\begin{figure}[hbt]
\centering
\subfigure[Urban area vs.\ sea and forest\label{ISF1}]{\includegraphics[width=.7\linewidth]{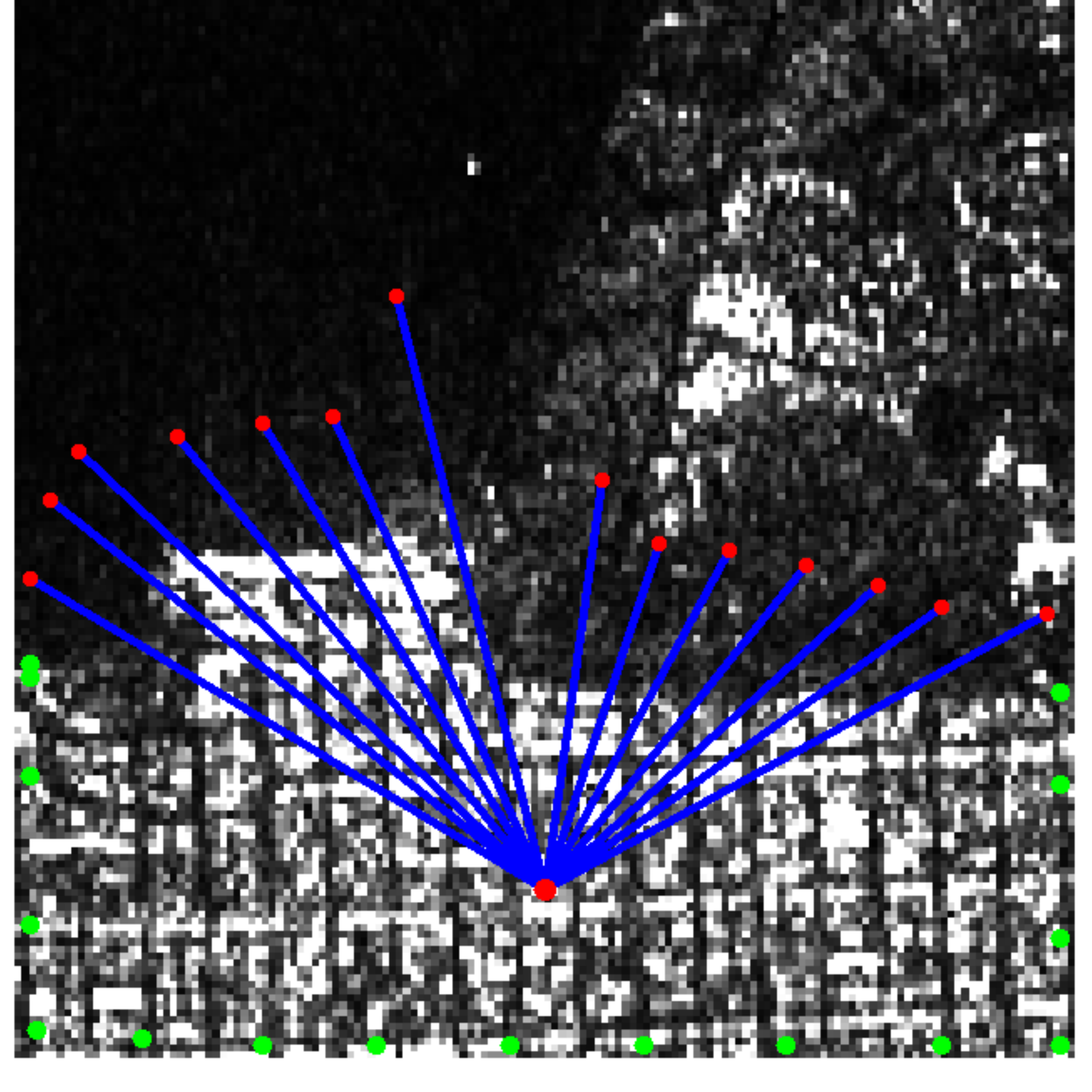}}
\subfigure[Detected egdes\label{ISF2}]{\includegraphics[width=1\linewidth]{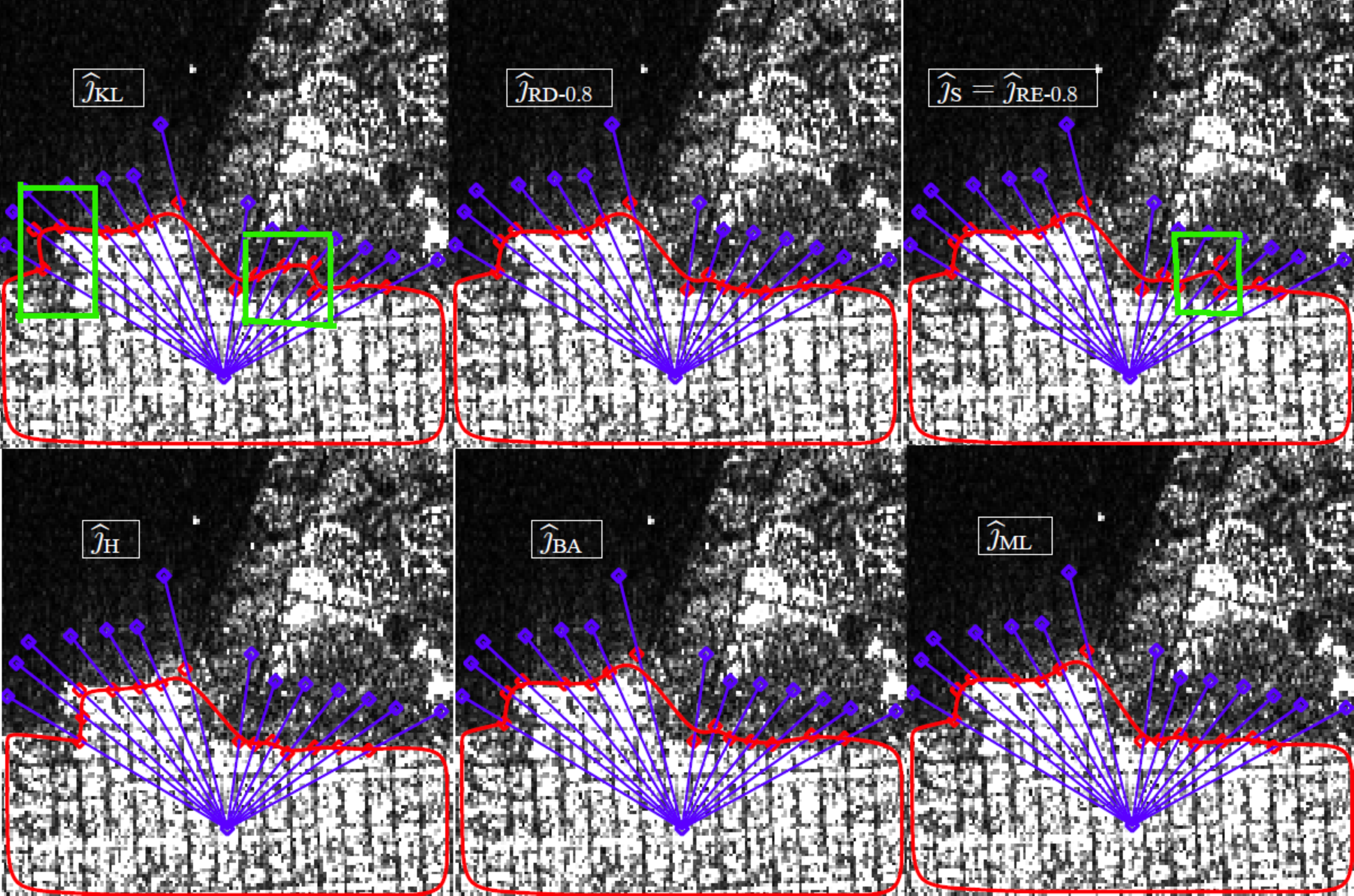}}
\caption{Performance of detectors with centroid on urban region from AIRSAR image of San Francisco Bay.
Errors highlighted in green boxes.
} 
\label{SimulationdiagonalII}
\end{figure}

\begin{figure}[hbt]
\centering
\subfigure[Sea vs.\ urban area and forest\label{IISF1}]{\includegraphics[width=.7\linewidth]{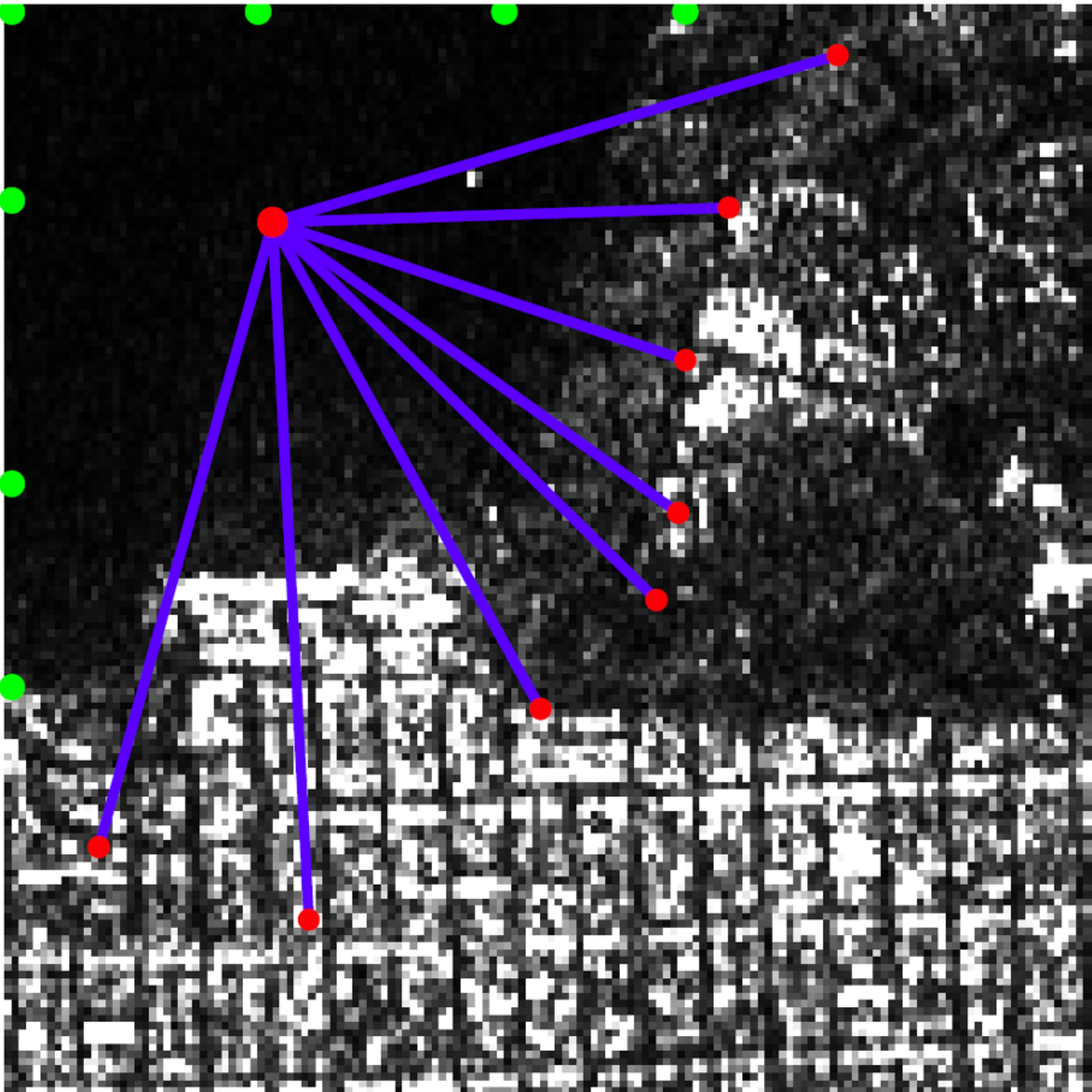}}
\subfigure[Detected egdes\label{IISF2}]{\includegraphics[width=1\linewidth]{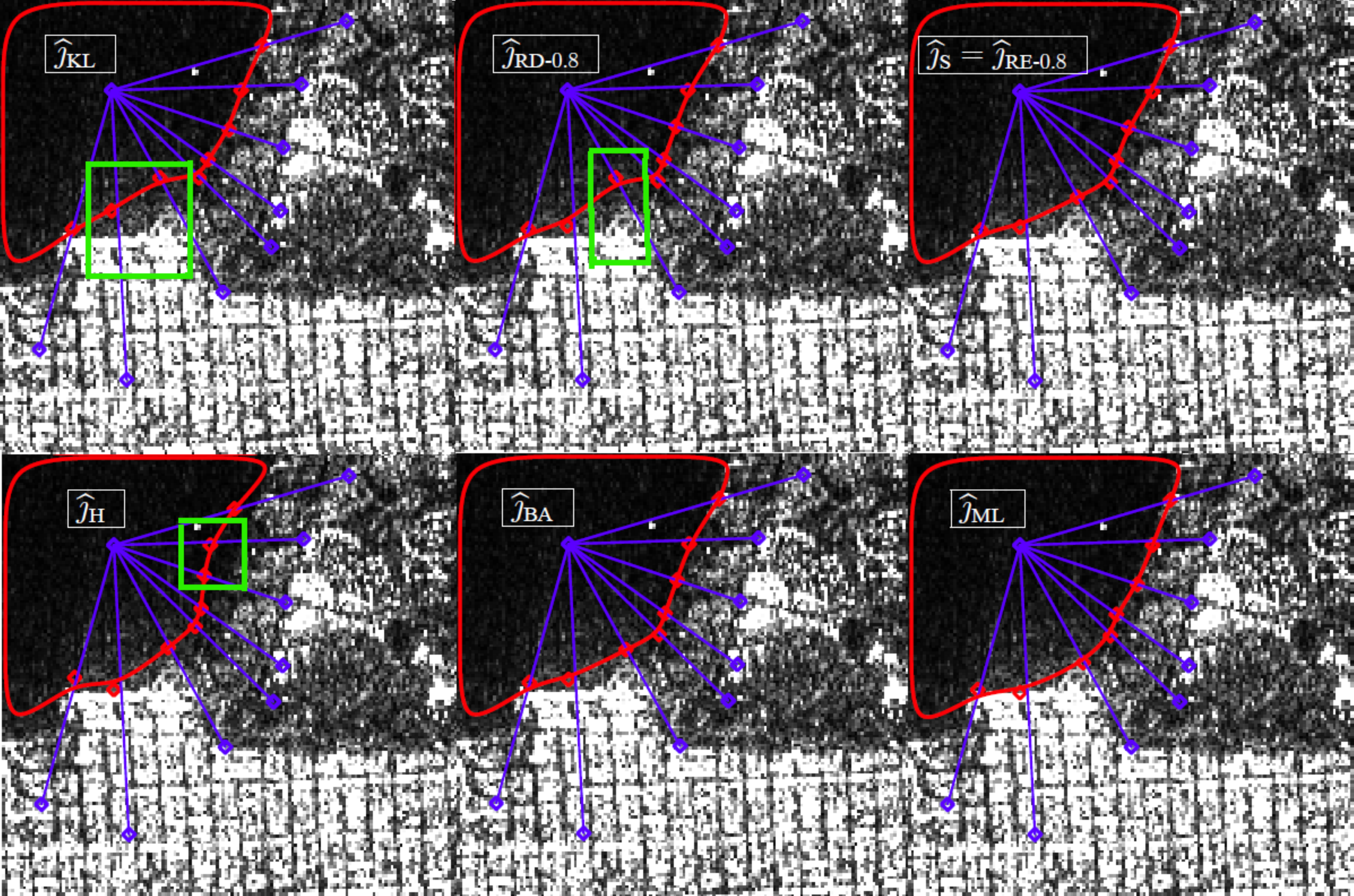}}
\caption{Performance of detectors with centroid on ocean region from AIRSAR image of San Francisco Bay.
Errors highlighted in green boxes.
} 
\label{SimulationdiagonalIII}
\end{figure}

The first application was performed on a San Francisco Bay image obtained by the AIRSAR sensor in L-band with four nominal looks, which is available in~\cite{PottierFACHPMWMLD09}.  
Fig.~\ref{ISF1} shows a $150\times150$ scene of the image with three well defined regions: ocean (dark area -- top and left part of the image), forest (gray area -- top and right), and urban area (light area -- bottom).
Figs.~\ref{SimulationdiagonalII} and~\ref{SimulationdiagonalIII} present the result of detecting the edges which separate urban and ocean regions from other classes, respectively.
The following analysis was performed by visual inspection on edges reconstructed from the estimated transition points and B-splines of fourth degree.

In Fig.~\ref{ISF2}, 
estimator $\widehat{\jmath}_{\text{KL}}$ presents two mildly biased boundaries (highlighted with green squares), whereas $\widehat{\jmath}_\text{S}$ and $\widehat{\jmath}_{\text{RE-}0.8}$ present only one. 
Remaining methods present similar behavior.

Fig.~\ref{IISF2} also exhibits three situations where the detected edge was slightly biased, namely $\widehat{\jmath}_{\text{KL}}$, $\widehat{\jmath}_{\text{RE-}0.8}$ and $\widehat{\jmath}_\text{H}$, each with one noticeable deviation.
Again, 
other discussed methods behave alike.

The Bhattacharyya and likelihood techniques did not miss any detection in these examples. 
This, and the simulation results on accuracy and computational time, leads us to suggest $\widehat{\jmath}_{\text{B}}$ as the fastest and most efficient detector.

\begin{figure}[hbt]
\centering
\subfigure[HH channel of the ESAR image with selected rays~\label{actualBD1}]{\includegraphics[width=1\linewidth]{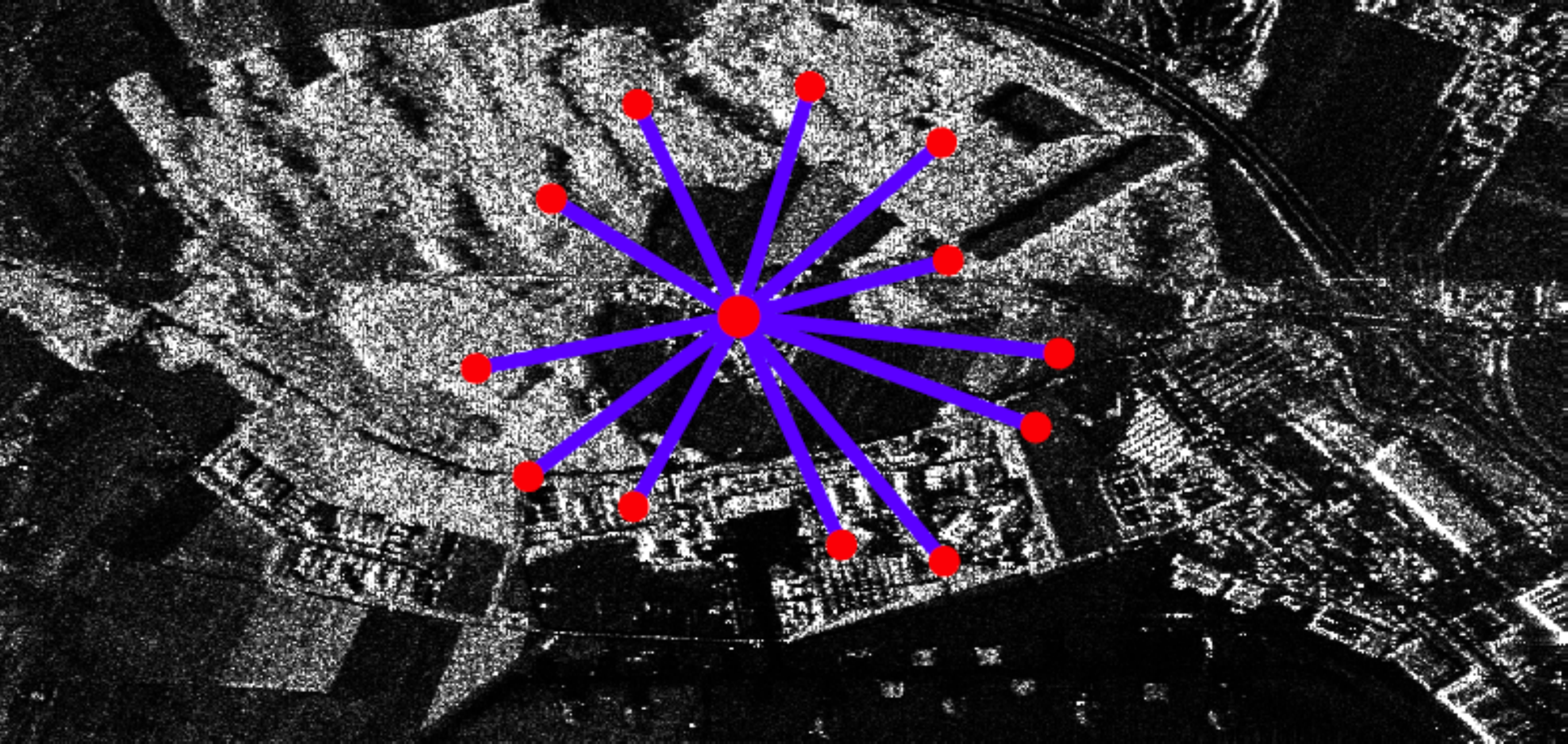}}\\
\subfigure[Detection in actual images \label{actualBDRV}]{\includegraphics[width=1\linewidth]{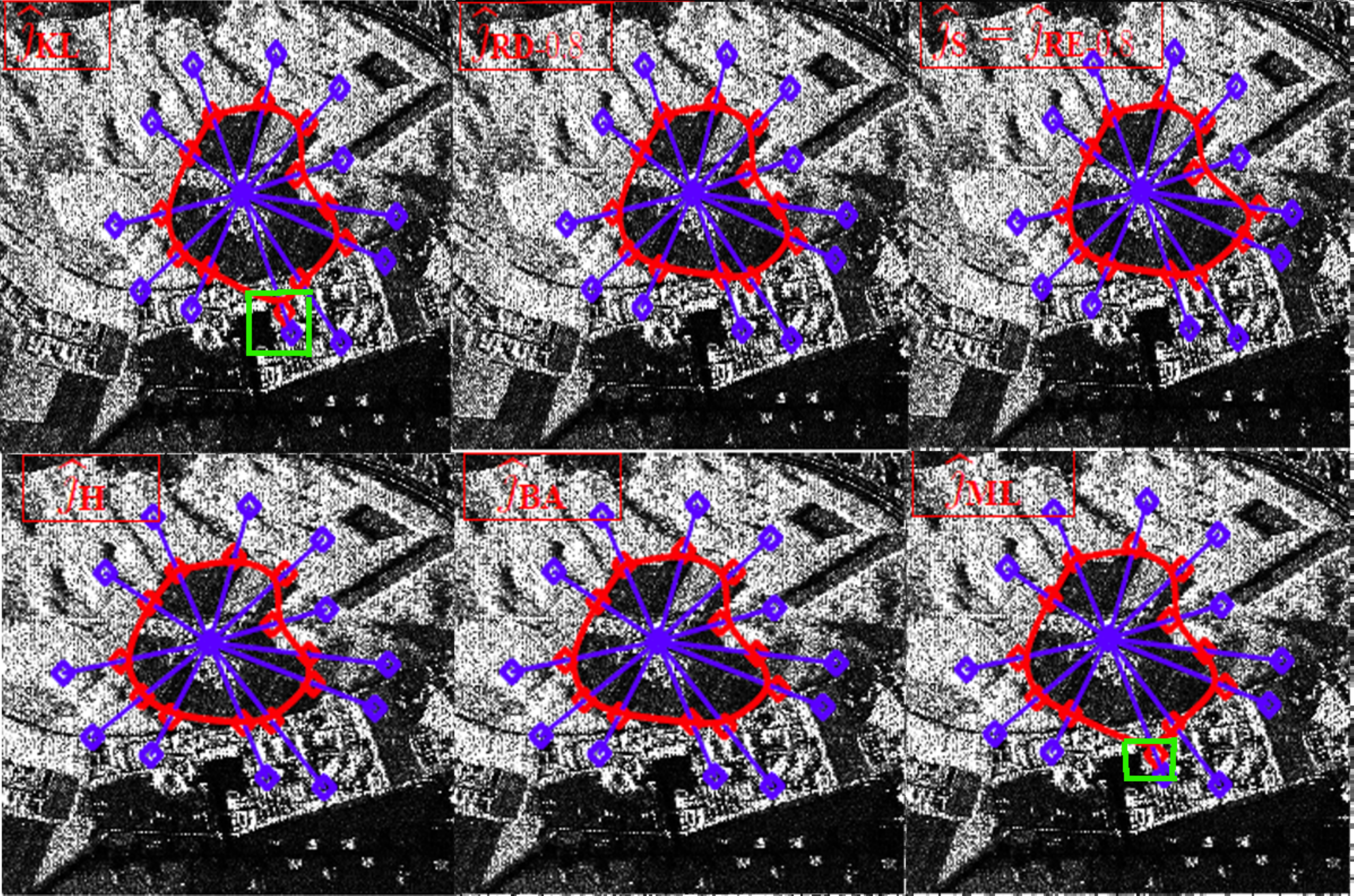}} 
\caption{Edge detection in an E-SAR image over the surroundings of Munich, Germany.
 Error highlighted in the green box.
} 
\label{SimulationdiagonalIV}
\end{figure}

Fig.~\ref{SimulationdiagonalIV} shows an E-SAR image over the surroundings of Munich, Germany, which  was obtained with $3.2$ (equivalent number of) looks.
The nonconvex region to the center of Fig.~\ref{actualBD1} is of interest, for which an arbitrarily centroid and thirteen rays were defined.
Fig.~\ref{actualBDRV} summarizes the results.
It shows the edges, and the label of the ray for which the transition point was identified far from the visual true point.
Both maximum likelihood and Kullback-Leibler methods failed at detecting only one situation, highlighted with a green square.
The other results are similar.

\begin{figure}[htbp]
\centering
\subfigure[Real Images with selected axes \label{actualBD2}]{\includegraphics[width=.48\linewidth]{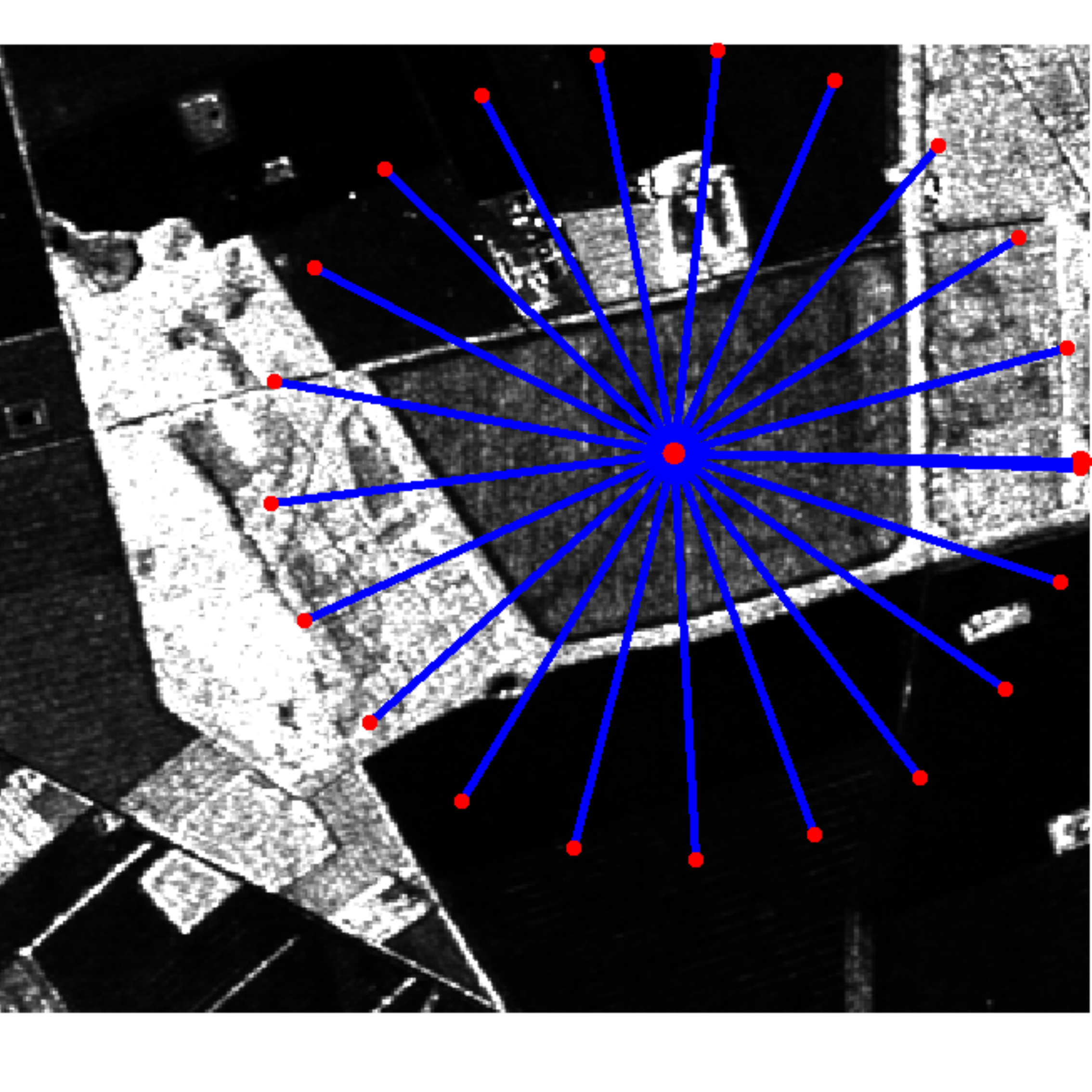}}
\subfigure[$\widehat{\jmath}_{\text{ML}}$ \label{actualBDRV2}]{\includegraphics[width=.48\linewidth]{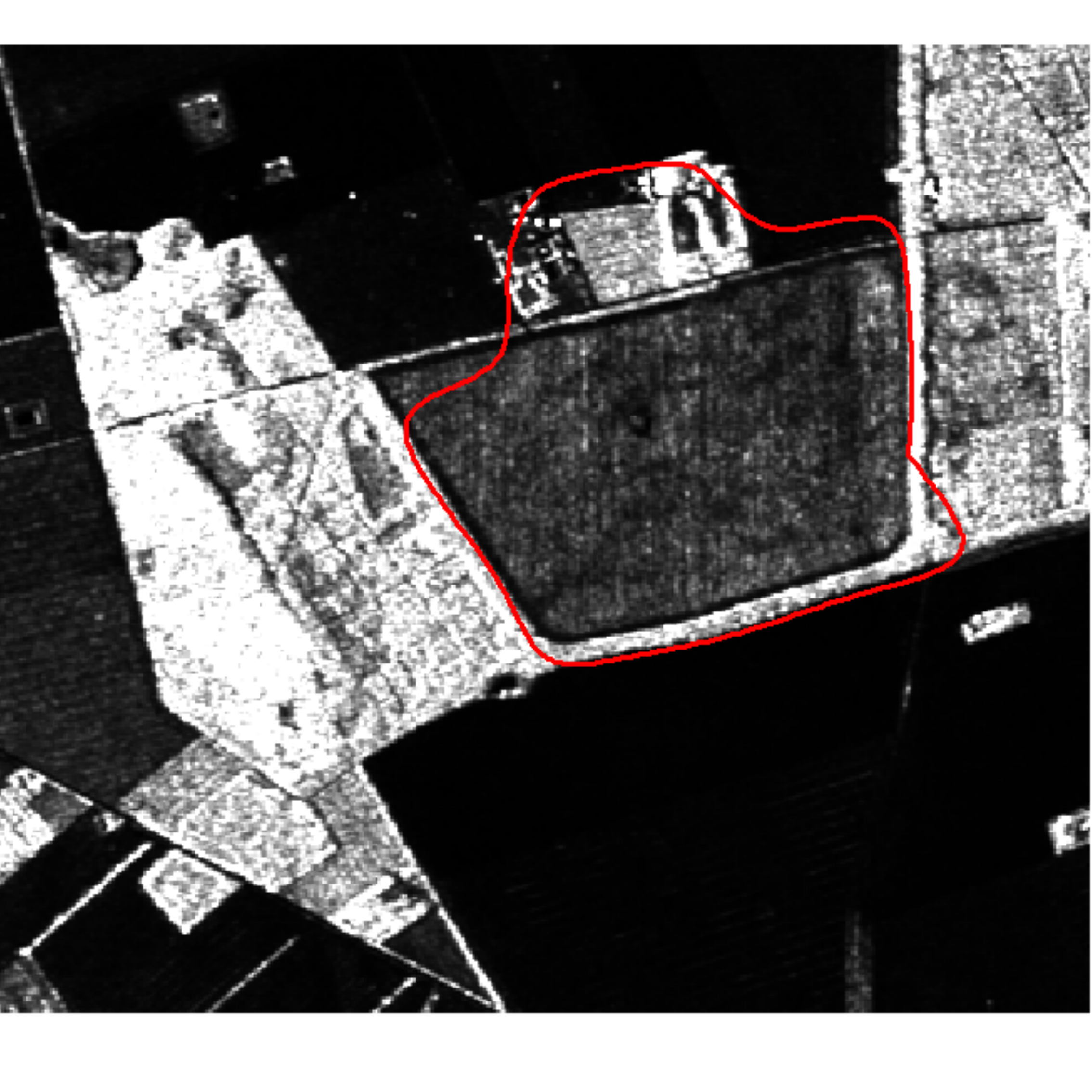}}
\subfigure[$\widehat{\jmath}_{\text{BA}}=\widehat{\jmath}_{\text{RD-}0.8}$ \label{actualBDHEL2}]{\includegraphics[width=.48\linewidth]{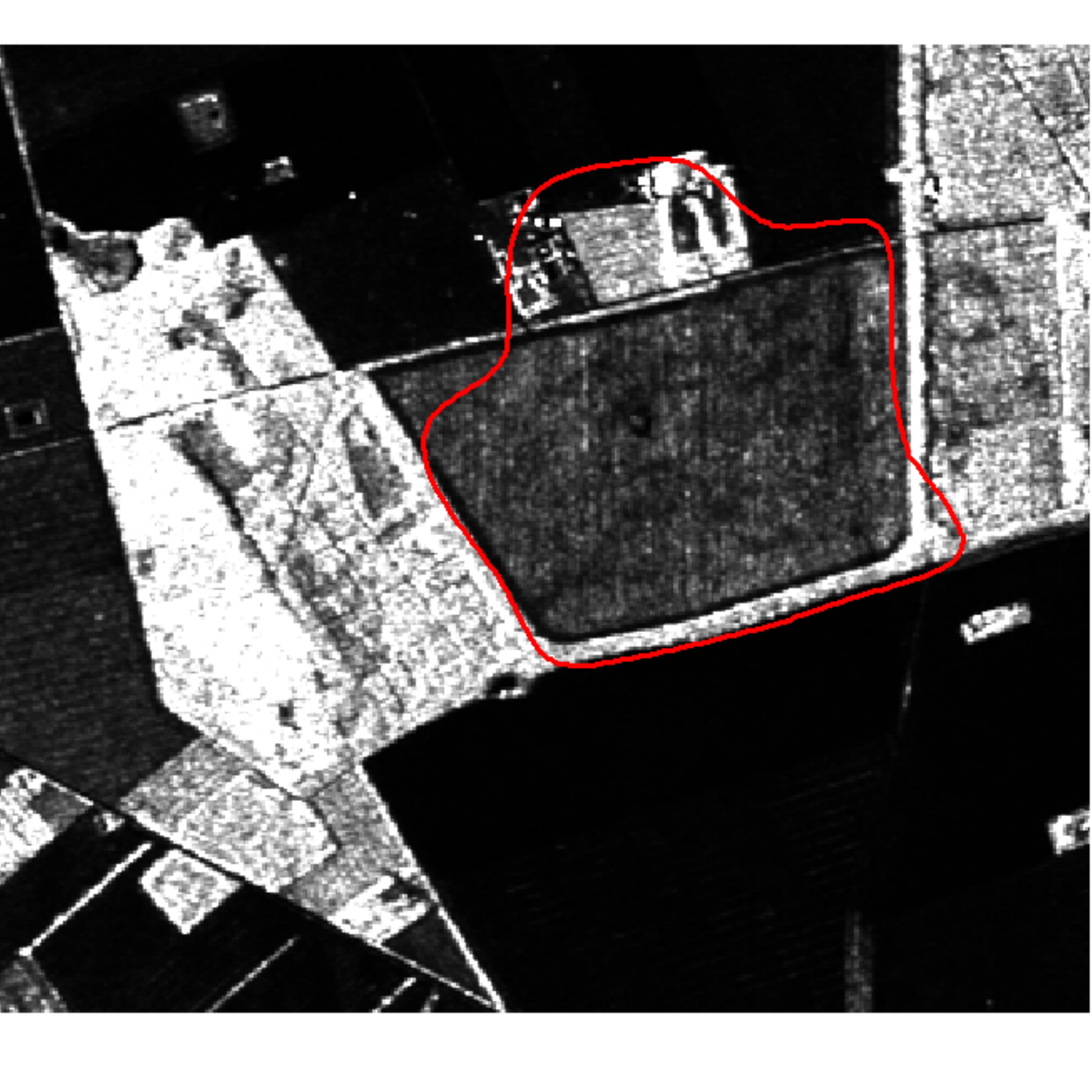}}
\subfigure[$\widehat{\jmath}_{\text{RE-}0.8}=\widehat{\jmath}_{\text{S}}$ \label{actualBDRV2DeNovo}]{\includegraphics[width=.48\linewidth]{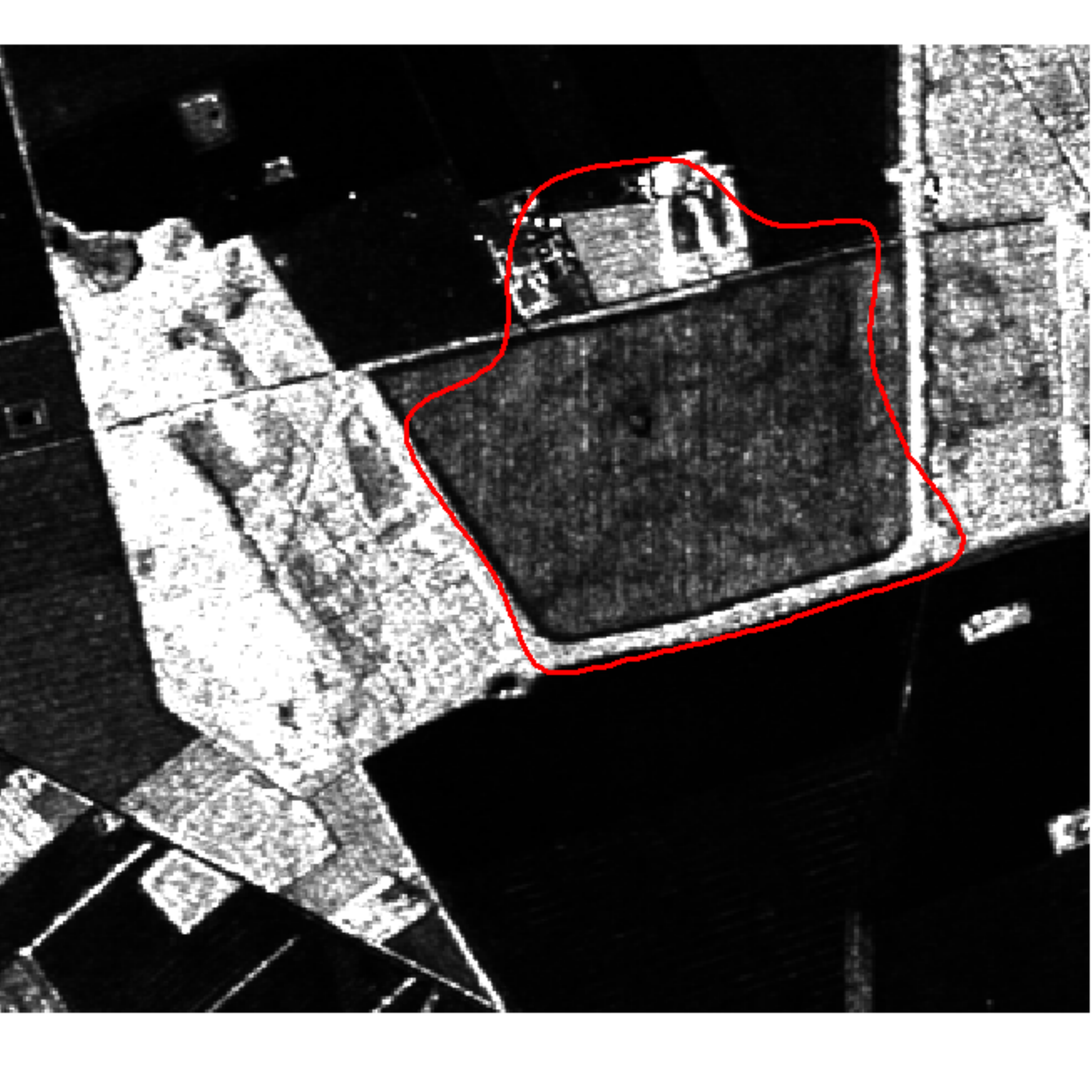}}
\subfigure[$\widehat{\jmath}_{\text{KL}}$ \label{actualBDKL2}]{\includegraphics[width=.48\linewidth]{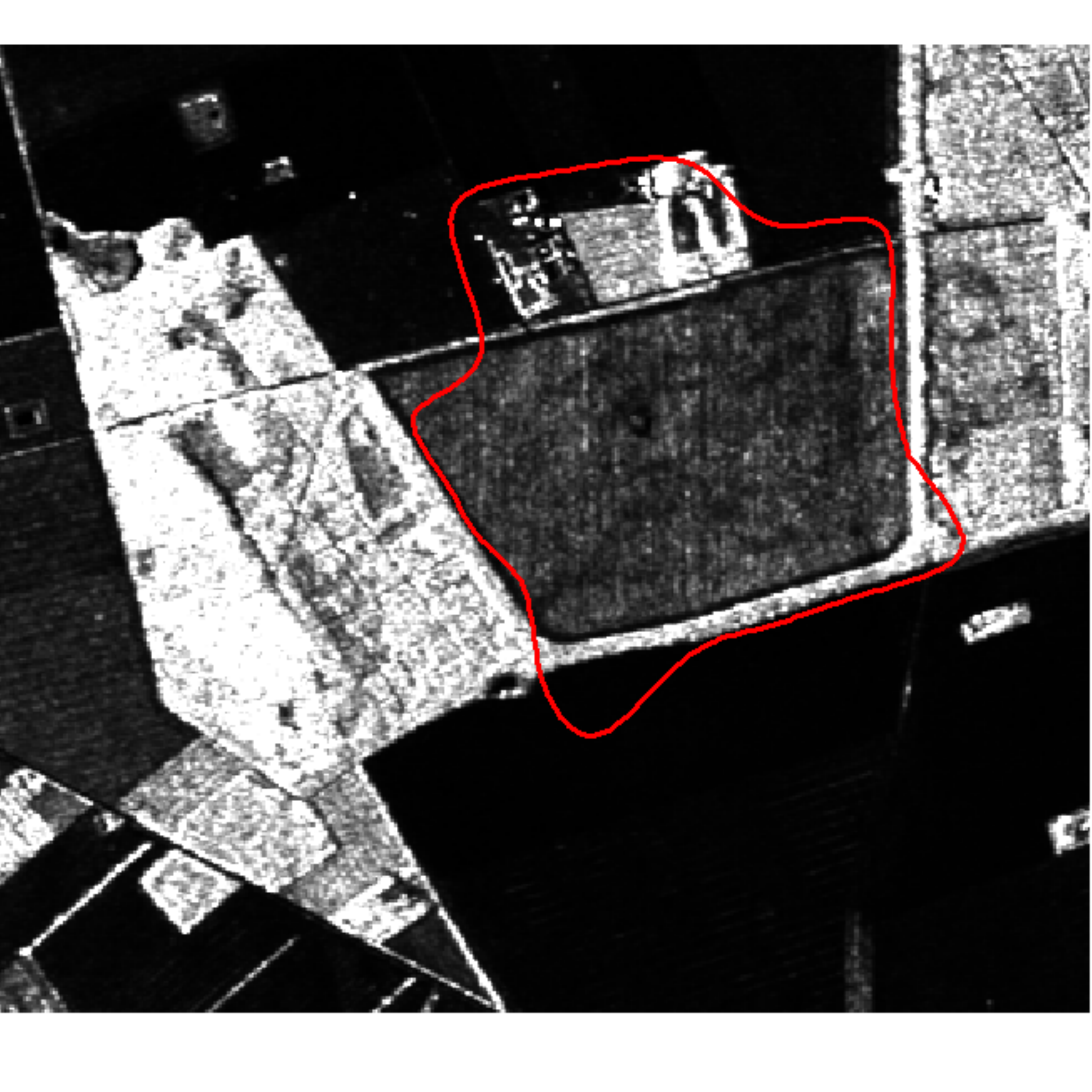}}
\subfigure[$\widehat{\jmath}_{\text{H}}$ \label{actualBDHEL3}]{\includegraphics[width=.48\linewidth]{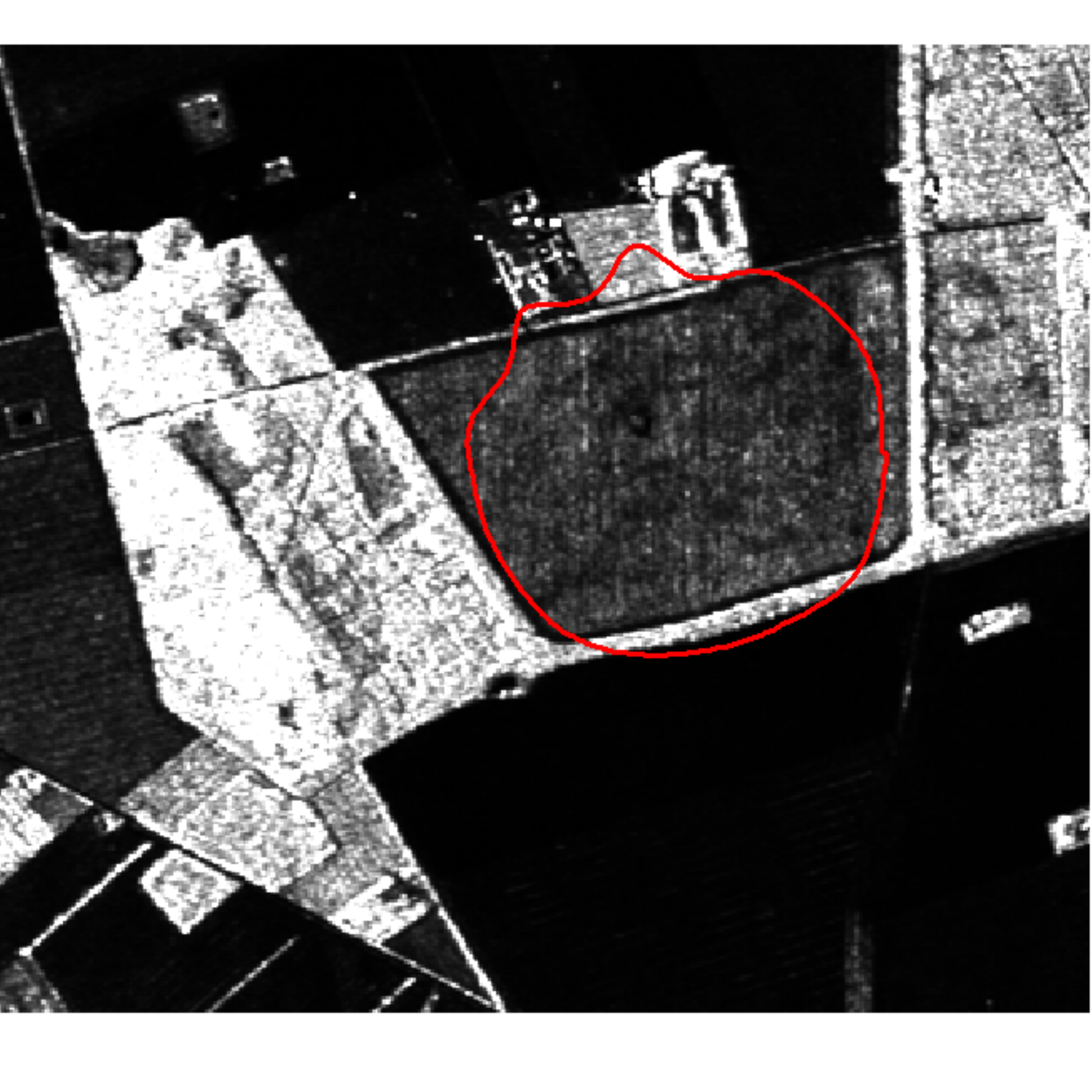}}
\caption{
Edge detection of an agricultural parcel in an EMISAR image
of Foulum.
} 
\label{actualpolsarimages2}  
\end{figure}

There are practical situations in which all methods perform alike.
Consider the EMISAR image of the agricultural area of Foulum, Denmark, presented in Fig.~\ref{actualpolsarimages2}.
This image was obtained in L-band and quad-pol; its HH channel is shown in Fig.~\ref{actualBD2}.
Notice that the method based on the joint likelihood (Fig.~\ref{actualBDRV2}) works similarly to $\widehat{\jmath}_{\text{BA}}$ and $\widehat{\jmath}_{\text{RD-}0.8}$ (Fig.~\ref{actualBDHEL2}),
$\widehat{\jmath}_{\text{RE-}0.8}$ (Fig.~\ref{actualBDRV2DeNovo}), 
and $\widehat{\jmath}_{\text{KL}}$ (Fig~\ref{actualBDKL2}).
However, the best technique regarding execution times is the one based on the Bhattacharyya distance.

\begin{figure}[htbp]
\centering
\subfigure[A ray crossing three different regions~\label{otherana1}]{\includegraphics[width=.76\linewidth]{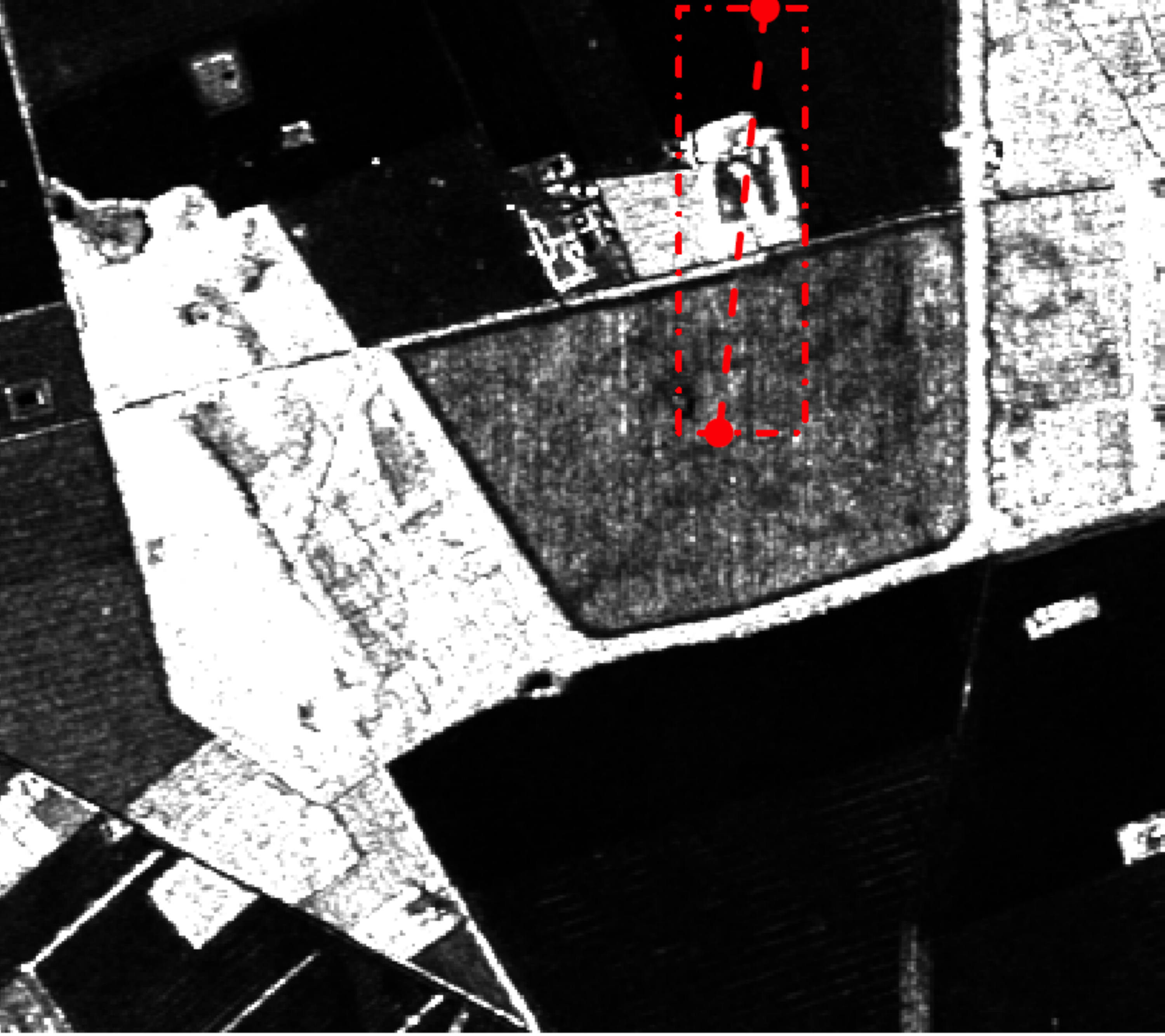}}
\subfigure[Zoom \label{otherana2}]{\includegraphics[width=.2\linewidth]{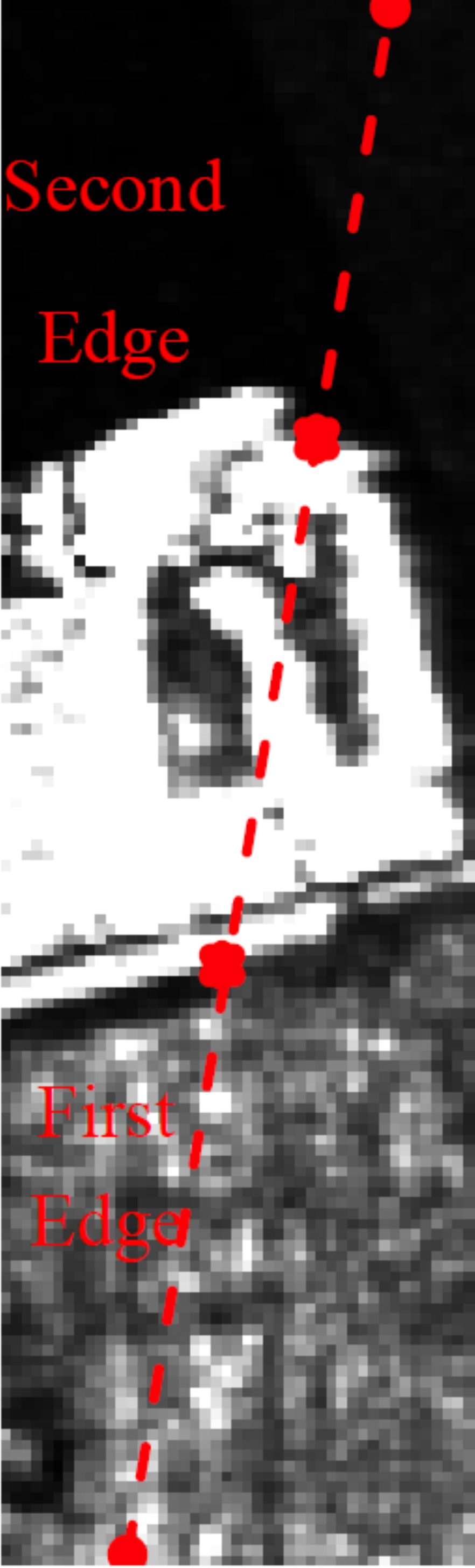}}
\subfigure[ $\widehat{\jmath}_{\text{ML}}$ \label{otherana3}]{\includegraphics[width=.48\linewidth]{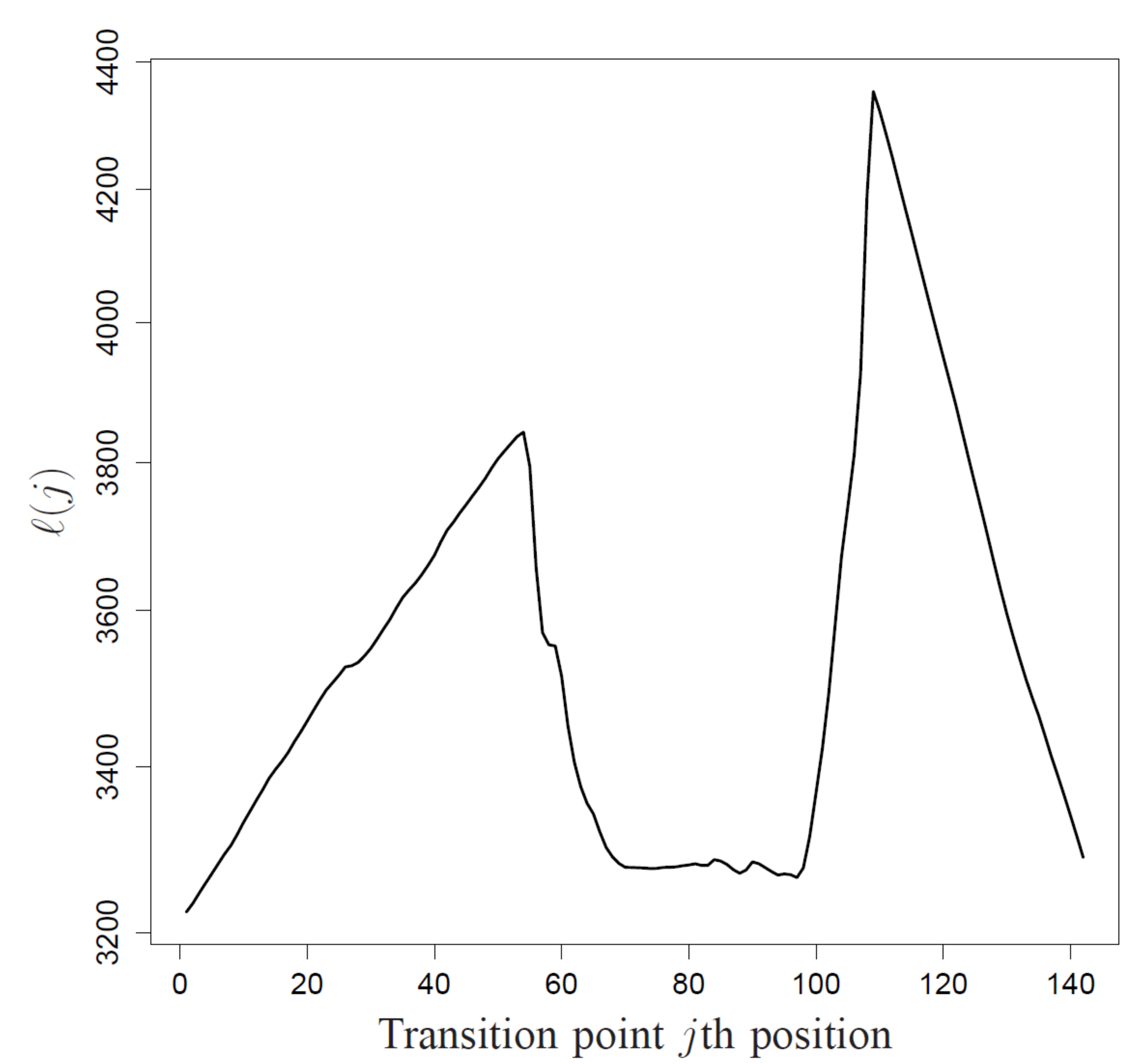}}
\subfigure[ $\widehat{\jmath}_{\text{H}}$ \label{otherana4}]{\includegraphics[width=.48\linewidth]{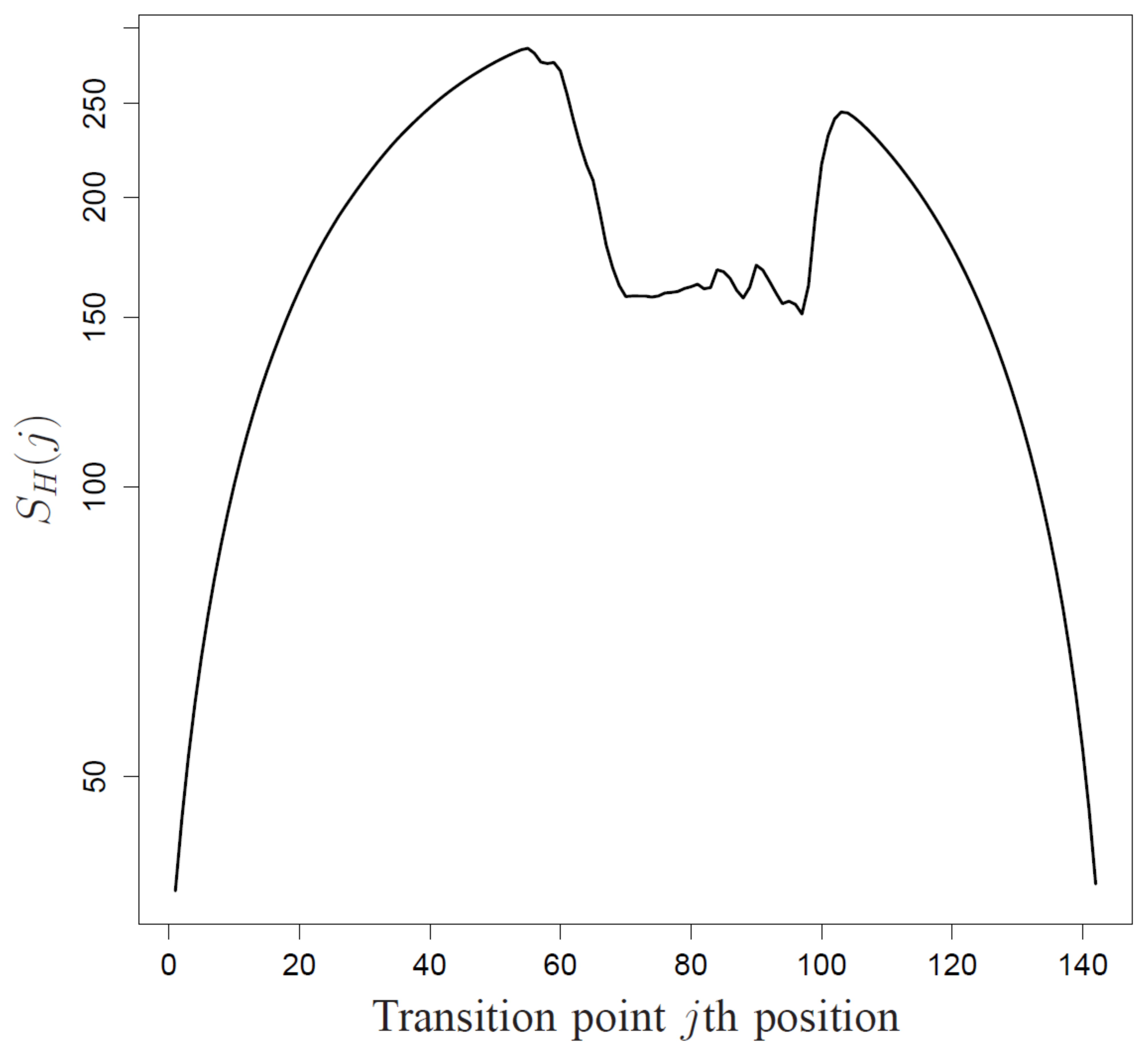}}
\subfigure[ $\widehat{\jmath}_{\bullet}$ \label{otherana5}]{\includegraphics[width=.48\linewidth]{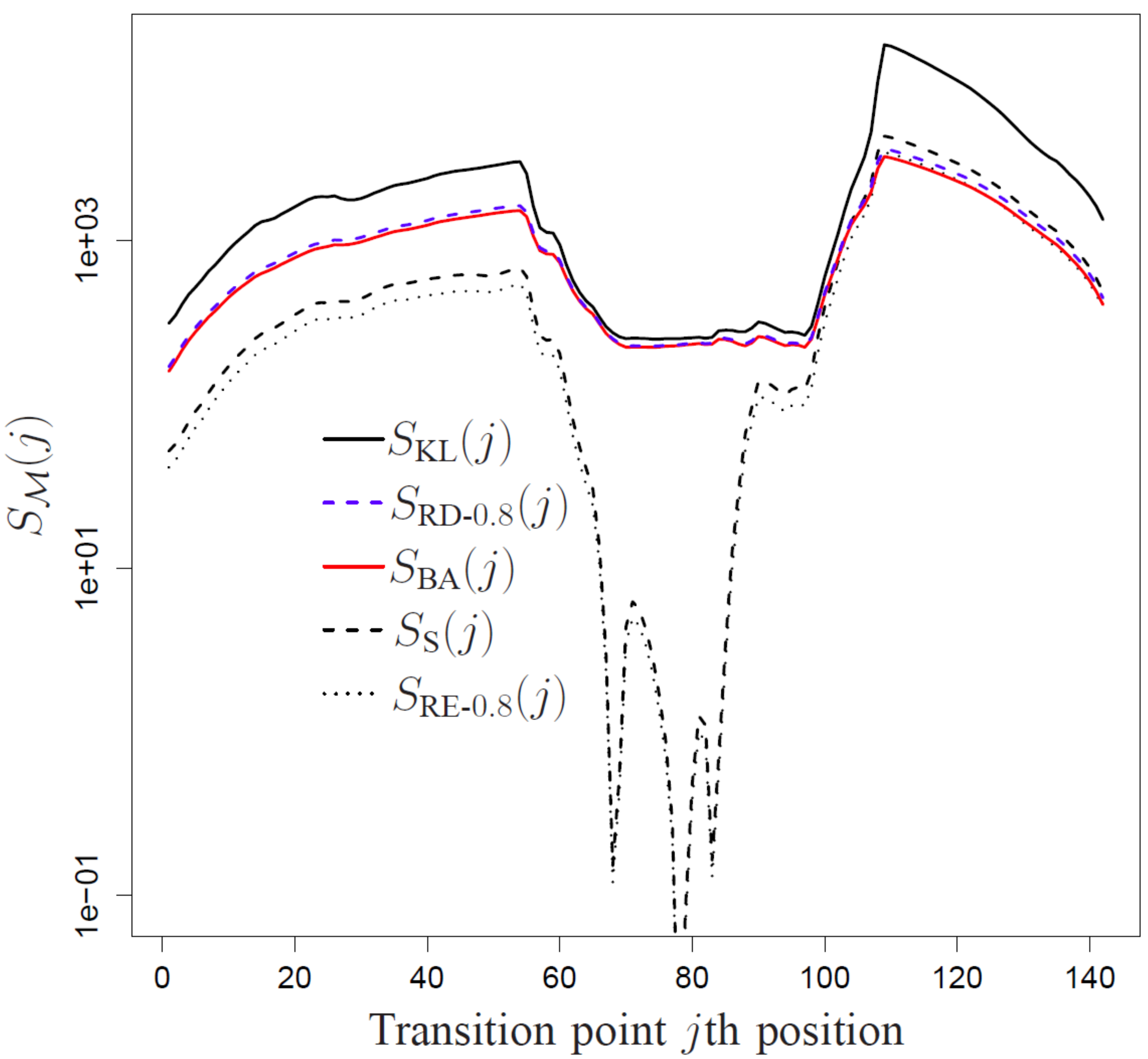}}
\caption{
Details of of the edge detection in a hard-to-deal-with
configuration from the EMISAR image of Foulum.
}
\label{otherana}
\end{figure}

In order to analyze further the result obtained by the Hellinger distance, presented in Fig.~\ref{actualBDHEL3}, consider the ray zoomed in Fig.~\ref{otherana2}, and notice that it crosses three regions: gray, white and black leading, therefore, to two possible edges.
As all the techniques are proposed in terms of the maximization of a discrimination measure, it is expected that the estimated boundary is the position of the highest distinction between two regions.
In this case, it takes place in the second edge, i.e., between the white and black areas.
This is confirmed by  Fig.~\ref{otherana3}, which shows the likelihood, but not in Fig.~\ref{otherana4}, which shows the Hellinger distance.
All other criteria behave alike, as shown in Fig.~\ref{otherana5}.

The Hellinger distance is limited to the $[0,1]$ interval. 
Additionally, let $v(\boldsymbol{k})=d_\text{H}(\boldsymbol{\theta}_1,\boldsymbol{\theta}_2 + \boldsymbol{k})$, empirical results provide evidence that there is a subset $\boldsymbol{\Theta}$ of the parametric space of the scaled complex Wishart distribution for which  $v(\boldsymbol{k})\approx 1$ for all $\boldsymbol{k}\in\boldsymbol{\Theta}$. 
This fact renders the statistics $S_\text{H}(\boldsymbol{\theta}_1,\boldsymbol{\theta}_2 + \boldsymbol{k})$ inappropriate for quantifying distinctions between distributed Wishart random fields when $\boldsymbol{k}\in \boldsymbol{\Theta}$.

This behavior opens new research lines based on a redefinition of the proposed methods in order to cope with situations in which rays cross three or more distinct regions.

\section{Conclusions}
\label{conclusion}

In this paper, boundary detection procedures using information theory measures (stochastic distances and entropies) and a likelihood function were considered. 
Due to its analytic tractability, these techniques are derived under the assumption that the Wishart distribution is an acceptable model for full PolSAR data.
In order to quantify and compare their performances, phantoms and real PolSAR images were used.

The performance of each procedure was quantified by the probability of correctly detecting the edge position within at most $k$ pixels, as well as by the execution time.
Problems with different levels of complexity were assessed in order to quantify the influence of the number of looks, the backscatter matrix as a whole, the SPAN, the covariance,
and the spatial resolution.
Also, the information provided by the full polarimetric format with respect to each intensity channel was assessed, leading to the conclusion that the complete data set is recommended for edge detection.
The results provide evidence that all methods are able to detect edges correctly within four pixels, with the exception of entropy-based techniques which fail in more complex scenarios with similar covariance matrices.

The discriminatory function based on the Kullback-Leibler distance is the fastest to compute, while the likelihood function is the slowest one among the considered measures.

The simulation studies revealed that information theory measures are consistently closer to the true boundary than the joint likelihood detector.       
In particular, Shannon and Kullback-Leibler detectors provided the most accurate detection and the smallest execution time, respectively.

Finally, three applications to real data obtained by E-SAR, EMISAR, and AIRSAR sensors were performed.
The seemingly poor performance of the edge detection based on the Hellinger distance was identified and explained by its properties.

As a conclusion, techniques which employ R\'enyi and Bhattacharyya distances and entropies outperform other methods and are recommended for edge detection in PolSAR imagery.

Further research aims at considering models which include heterogeneity~\cite{FreitasFreryCorreia:Environmetrics:03,FreryCorreiaFreitas:ClassifMultifrequency:IEEE:2007,
PolarimetricSegmentationBSplinesMSSP}, robust, improved and nonparametric inference~\cite{AllendeFreryetal:JSCS:05,NonparametricEdgeDetectionSpeckledImagery,SilvaCribariFrery:ImprovedLikelihood:Environmetrics,
VasconcellosFrerySilva:CompStat}, and small samples techniques~\cite{FreryCribariSouza:JASP:04}.

\section*{Acknowledgment}

The authors are grateful to CNPq, Facepe, Fapeal, FAPESP, and Capes for their support.

\bibliographystyle{IEEEtran}
\bibliography{JSTARS2012}

\begin{IEEEbiography}[{\includegraphics[width=1in]{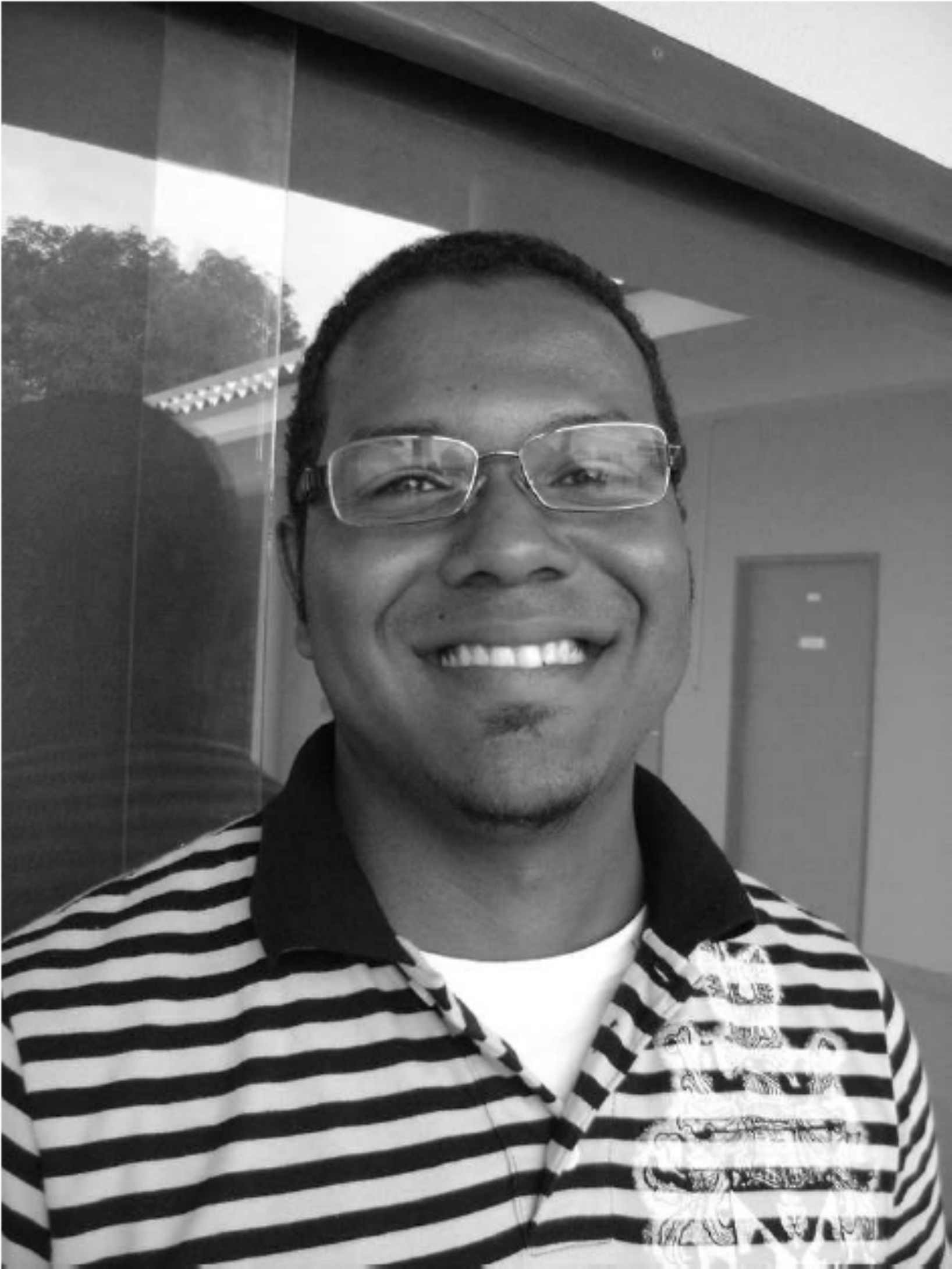}}]{Abra\~ao D.\ C.\ Nascimento}
holds B.Sc.\, M.Sc.\, and D.Sc. degrees in Statistics from Universidade Federal de Pernambuco (UFPE), Brazil, in 2005, 2007, and 2012, respectively.
In 2013, he joined the Department of Statistics at Universidade Federal da Para\'iba.
His research interests are statistical information theory, inference on random matrices, and asymptotic theory.
\end{IEEEbiography}

\begin{IEEEbiography}[{\includegraphics[width=1in]{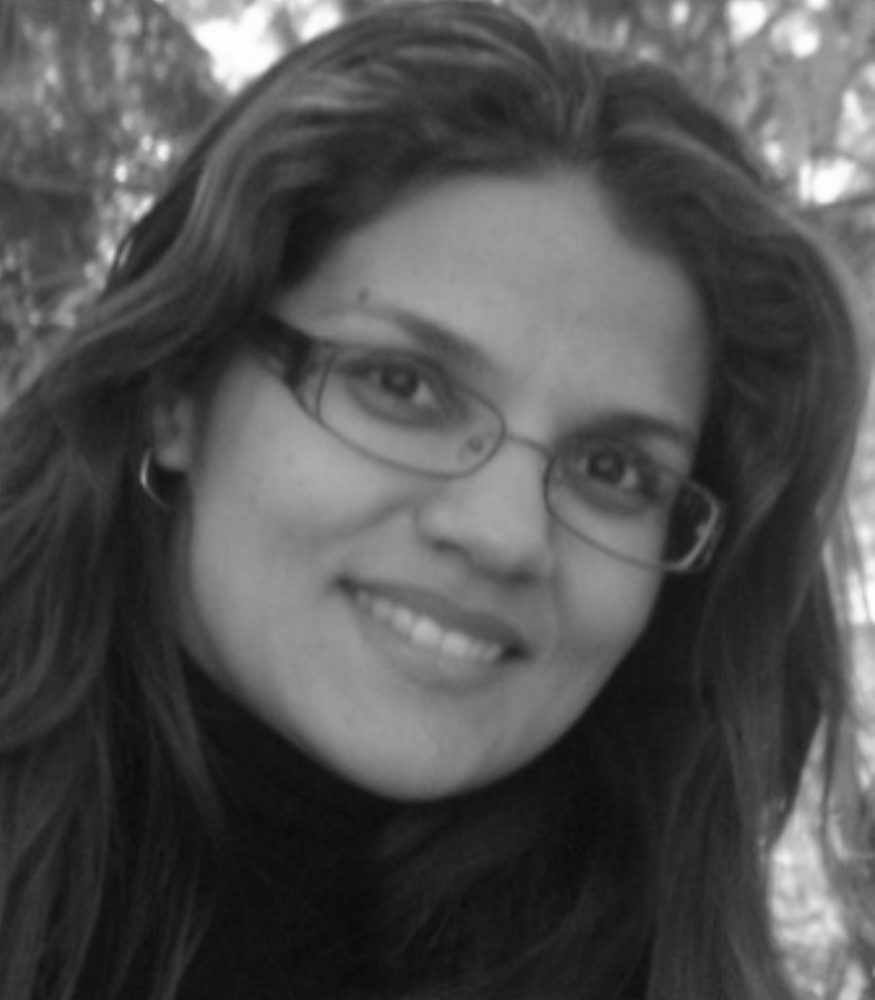}}]{Michelle M.\ Horta}
received the B.S. degree in computer sciences from Universidade Cat\'olica de Pernambuco, Recife, Brazil, in 2002, the M.Sc. degree in computer sciences from Universidade Federal de Pernambuco, Recife, Brazil, in 2004, and the Ph.D. degree in applied physics from Universidade de S\~ao Paulo, S\~ao Carlos, Brazil, in 2009. She is currently working toward the Postdoctoral Research in the Departamento de Computa\c{c}\~ao, Universidade Federal de S\~ao Carlos, S\~ao Carlos, Brazil. Her research interests are remote sensing, pattern recognition and statistical models.
\end{IEEEbiography}

\begin{IEEEbiography}[{\includegraphics[width=1in]{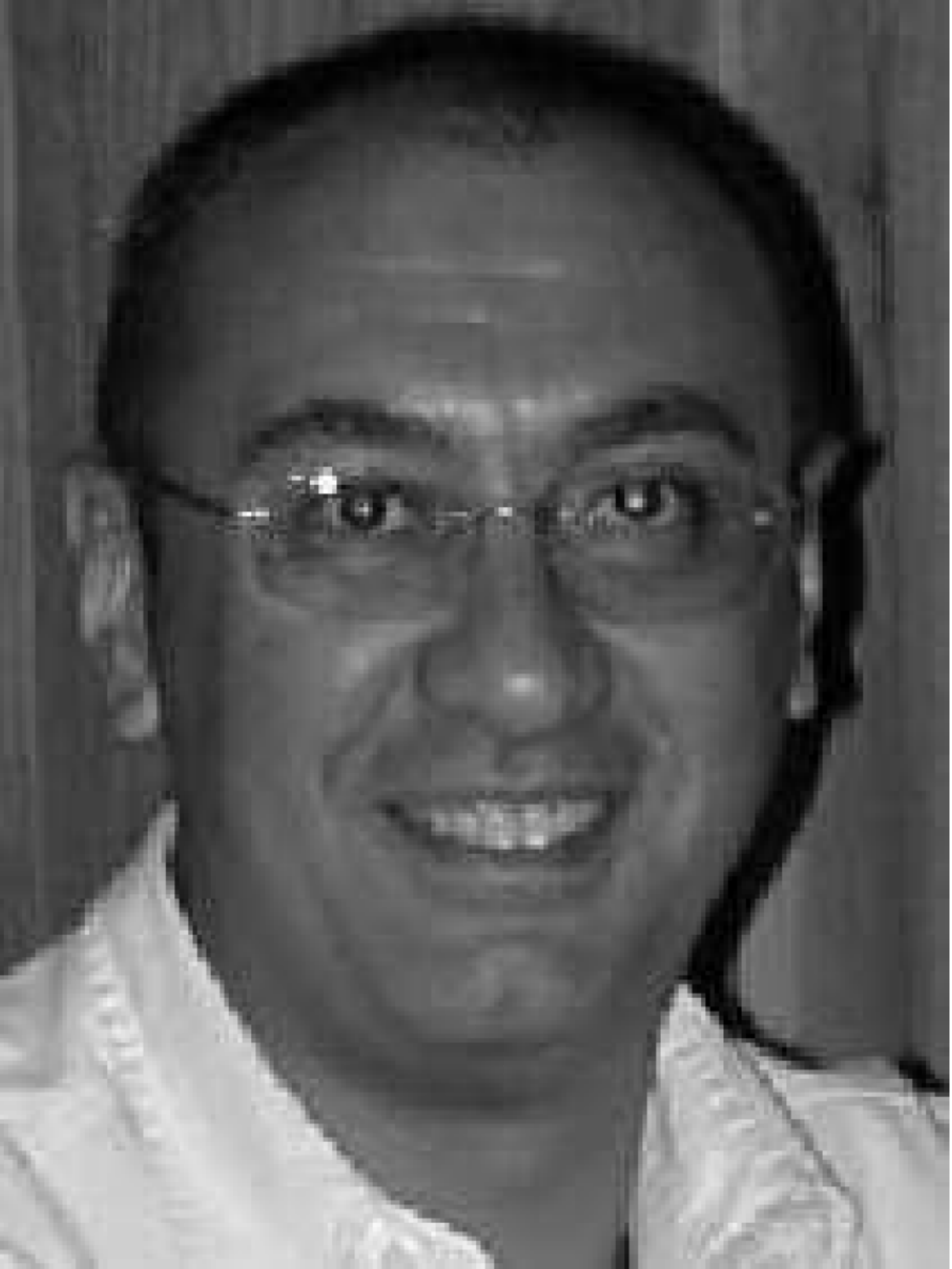}}]{Alejandro C.\ Frery}
graduated in Electronic and Electrical Engineering from the Universidad de Mendoza, Argentina.
His M.Sc. degree was in Applied Mathematics (Statistics) from the Instituto de Matem\'atica Pura e Aplicada (Rio de Janeiro) and his Ph.D. degree was in Applied Computing from the Instituto Nacional de Pesquisas Espaciais (S\~ao Jos\'e dos Campos, Brazil).
He is currently with the Instituto de Computa\c c\~ao, Universidade Federal de Alagoas, Macei\'o, Brazil.
His research interests are statistical computing and stochastic modeling.
\end{IEEEbiography}

\begin{IEEEbiography}[{\includegraphics[width=1in]{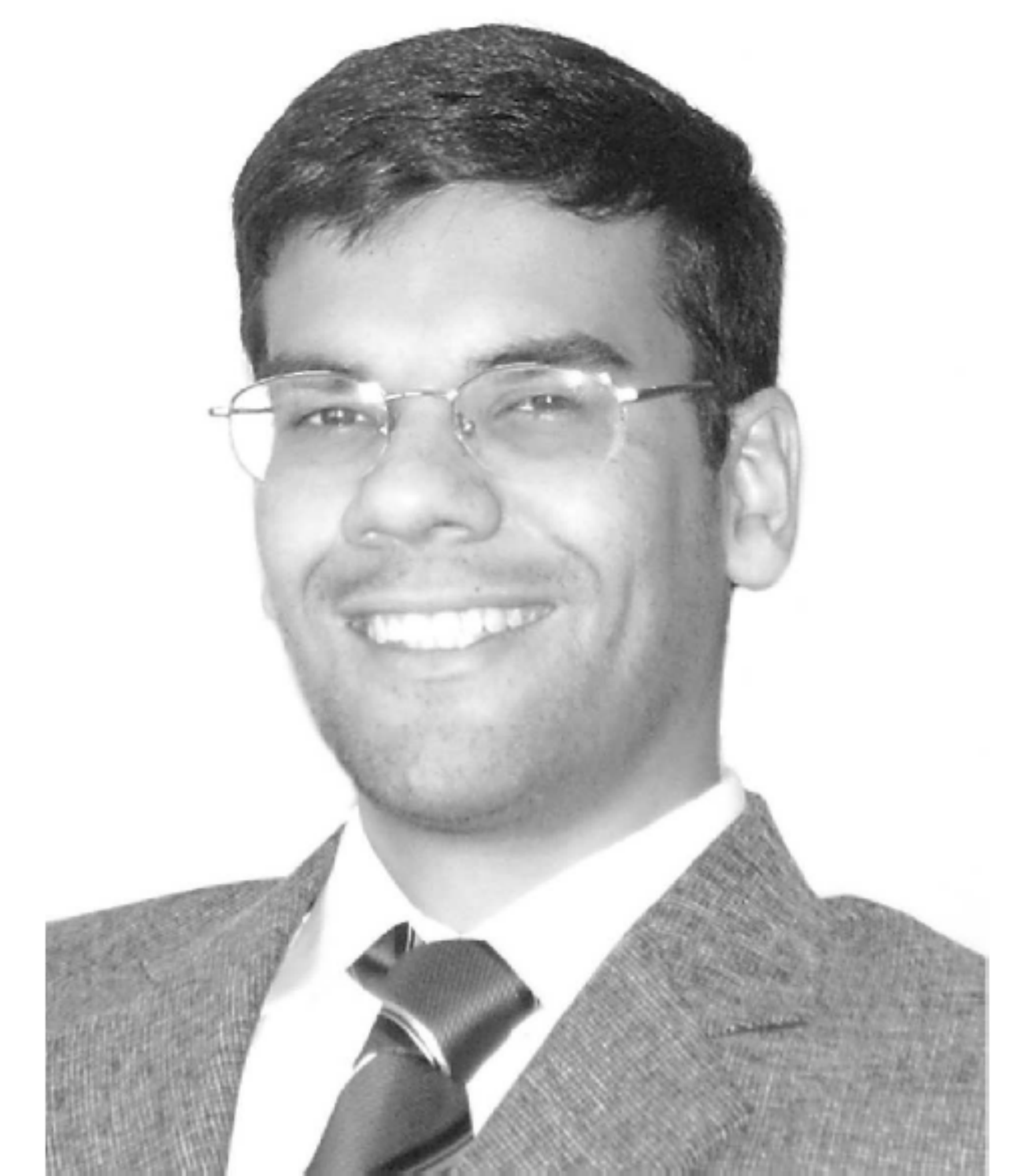}}]{Renato J.\ Cintra}
earned his B.Sc., M.Sc., and D.Sc. degrees in Electrical Engineering from
Universidade Federal de Pernambuco, Brazil, in 1999, 2001, and 2005, respectively.
In 2005, he joined the Department of Statistics at UFPE.
During 2008-2009, he worked at the University of Calgary, Canada, as a visiting research fellow.
He is also a graduate faculty member of the Department of Electrical and Computer Engineering, University of Akron, OH.
His long term topics of research include theory and methods for digital signal processing, communications systems, and applied mathematics.
\end{IEEEbiography}

\vfill

\end{document}